\pdfoutput=1
\documentclass[11pt,a4paper]{amsart}
\usepackage{graphics}
\usepackage{latexsym}
\usepackage{verbatim}
\usepackage{amsmath}
\usepackage{amsthm}
\usepackage{amssymb}
\usepackage{MnSymbol}
\usepackage{stmaryrd}

\usepackage[]{hyperref}
\hypersetup{
    colorlinks=true,%
    filecolor=black,%
    linkcolor=black,%
    urlcolor=black,  %
    citecolor = red
}
\usepackage [table]{xcolor}
\usepackage{multirow}
\usepackage{float}
\usepackage{tikz}
\usepackage{subfig}
\usepackage{algorithm2e}

\graphicspath{{pix/}}	
\newtheorem{theorem}{Theorem}

\newtheorem{proposition}[theorem]{Proposition}

\newcommand\CC{\hbox{C\kern -.58em {\raise .54ex \hbox{$\scriptscriptstyle |$}}
  \kern-.55em {\raise .53ex \hbox{$\scriptscriptstyle |$}} }}
\newcommand\NN{\hbox{I\kern-.2em\hbox{N}}}
\newcommand\RR{\mathbb{R}}
\newcommand\ZZ{{{\rm Z}\kern-.28em{\rm Z}}}

\newcommand\Div{ \textrm{div}}

\newcommand\xx{ \mathbf{x} }

\usepackage{algorithmic}

\begin{document}

\title[Mixed semi-Lagrangian/finite difference methods  for plasma simulations]{Mixed semi-Lagrangian/finite difference methods  for
  plasma simulations}

\author{ Francis Filbet}

\address{
Universit\'e de Lyon \& Inria,
Institut Camille Jordan, EPI Kaliffe,
43 boulevard 11 novembre 1918,
F-69622 Villeurbanne cedex,  FRANCE}

\email[F.~Filbet]{ filbet@math.univ-lyon1.fr }

\author{Chang Yang}

\address{Department of Mathematics, Harbin Institute of Technology, 92 West Dazhi Street, Nan Gang District,
Harbin 150001, China}

\email[C.~Yang]{yangchang@hit.edu.cn}


\maketitle

\begin{abstract}

In this paper, we present an efficient algorithm for the long time
behavior of  plasma simulations. We will focus on 4D drift-kinetic
model, where the plasma's motion occurs  in the plane perpendicular to
the magnetic field and can be governed by the 2D guiding-center
model. 

Hermite WENO reconstructions, already  proposed in \cite{YF15}, are applied for solving the
Vlasov equation.  Here we consider an arbitrary computational domain with an appropriate numerical method for the treatment of boundary conditions. 

Then we apply this algorithm for plasma turbulence simulations. We
first solve the 2D guiding-center model in a D-shape domain and
investigate the numerical stability of the steady state. Then,  the 4D
drift-kinetic model is studied  with a mixed method, {\it i.e.} the
semi-Lagrangian method in linear phase and finite difference method
during the nonlinear phase. Numerical results show that the mixed
method is efficient and accurate in linear phase and  it is much
stable during the nonlinear phase. Moreover, in practice  it has
better conservation  properties.

\end{abstract}

\vspace{0.1cm}


\vspace{0.1cm}

\noindent 
{\small\sc 2000 Mathematics Subject Classiﬁcation.}  {\small   65M08, 65M25, 78A35.}


\tableofcontents


\section{Introduction} 
\setcounter{equation}{0}
\label{sec:Intro}

  In the context of strongly magnetized plasma simulations, the motion
  of particles  is confined around the magnetic field lines;  the
  frequency of this cyclotron motion is faster than the frequencies of
  interest. Therefore, the physical system can be reduced from the
  $6D$ Vlasov-Maxwell system to a four or  five dimensional model by
  averaging over the gyroradius of charged particles (See for a review
  \cite{ref1, ref2}).  In this paper we focus on 4D drift-kinetic
  model, where the movement of the plasma in the plane perpendicular
  to the magnetic field can be governed by the guiding-center
  model. This reduced model could help us to investigate plasma
  turbulence problems with an acceptable computational time. More especially, using the 2D guiding-center model, we can focus directly on the difficulties of treatment of boundary conditions on arbitrary computational domain. Moreover, these reduced models have the conservative properties, which can be used as criterion to evaluate the good performance of numerical methods.

To develop accurate and stable numerical techniques for plasma
turbulence (4D drift kinetic, 5D gyrokinetic and 6D kinetic models) is
one of our objectives. In \cite{bibGV, ref5} several numerical solvers
have been developed using an Eulerian formulation for gyro-kinetic
models. However, spurious oscillations often appear  in the nonlinear
phase  when small structures occur and it is difficult to distinguish
physical and numerical oscillations. Moreover, for these models
semi-Lagrangian methods are no more conservative, hence the long time
behavior of the numerical solution may become unsuitable.  At
contrast,  a class of numerical methods based on the Hermite
interpolation \cite{bibFS}, together with a   weighted essentially
non-oscillatory (WENO)  reconstruction can be applied either to
semi-Lagrangian or to finite difference methods \cite{YF15}.

On the one hand, the semi-Lagrangian methods are very efficient and
fast but do not conserve  mass in an arbitrary grid. On the other
hand, the  finite difference methods are mass conservative and more
stable for long term simulations but have a restrictive  CFL constraint. We thus propose a mixed method to carry forward the advantages of each method, {\it i.e.} the semi-Lagrangian method in linear phase where the solution is relatively smooth, the  finite difference method during the nonlinear phase where a lot of small structures appear. We will apply the mixed method to the  4D drift-kinetic model to verify its efficiency.

 The numerical resolution of 4D drift-kinetic model in a cylinder has
 been already performed  via polar coordinates~\cite{bibGV}.  However,
 polar coordinates are not always suitable since artificial singular
 points appear  in the computational domain  coming from the change of
 variable. Moreover, for more complicated geometry, it is not
 straightforward to perform the appropriate  change of variables. From
 these considerations, we focus here   on the discretization of
 transport models on a Cartesian mesh and develop a suitable technique
 as in~\cite{bibFY1} to treat boundary conditions. This method  is
 based on the inverse Lax-Wendroff method \cite{bibFY1}. On the other
 hand, to compute the electric potential from the Poisson equation, we
 apply an extrapolation technique \cite{bibFY2}. To test the numerical
 algorithm, the guiding center model will be solved  on a $D$-shape
 domain~\cite{bibmiller}. We will first compute  a steady state
 solution by solving numerically a nonlinear Poisson equation  and
 then investigate its stability properties.

The paper is organized as follows : in  Section~\ref{sec:model}, we
present the derivation of the  4D drift-kinetic and the 2D
guiding-center models, and  their conservative properties.  Then in
Section \ref{sec:schemes},  we recall the Hermite WENO reconstructions
developed in~\cite{YF15} for solving the Vlasov equations, and  the treatment of boundary conditions
corresponding to  the Poisson equation. In Section \ref{sec:Num}, we
numerically compute a steady state solution for the 2D guiding-center
model in a $D$-shape domain. Then, we study its stability properties
by performing  numerical simulations with a perturbed steady state as
initial data. Finally,  we numerically solve the 4D drift-kinetic model with the mixed methods, and present the numerical results of the conservative properties of the 4D drift-kinetic model and the evolution of the distribution function of density.

\section{Mathematical models}
\label{sec:model}
\setcounter{equation}{0}

The Vlasov equation for the distribution function $f$ in standard form in
standard notation  is
\begin{equation}
  \frac{\partial f}{\partial t}+\mathbf{v}\cdot\nabla_{\mathbf{x}}f+\frac{e}{m}\left(\mathbf{E}+\frac{\mathbf{v}\times\mathbf{B}}{c}\right)=0,
  \label{eq:vlasov}
\end{equation}
where $t\in\mathbb{R}^+$ is the time variable,
$\mathbf{x}\in\Omega\subset\mathbb{R}^3$ is the space variable, 
$\mathbf{v}\in\mathbb{R}^3$ is the  velocity variable, $m$ is the particle
mass, $e$ is its charge, $\mathbf{E}$ is the electric field and
$\mathbf{B}$ is the magnetic field.
We assume the electric field   is computed by
$\mathbf{E}=-\nabla\phi$, where $\phi$ is electric potential whereas the magnetic field    is uniform $\mathbf{B}=B\,e_z$, where $ e_z$ stands for the unit vector in the toroidal direction.
Moreover, we assume that $f$ is vanishing at infinity of velocity field and periodic boundary condition is taken in $z$ direction.

To derive the drift kinetic model, we start to perform a change of variable according  to the drift direction $\mathbf{E}\times\mathbf{B}$; that is
\begin{equation}
 \mathbf{w}=\begin{pmatrix}
             \mathbf{w}_\bot\\\mathbf{w}_\|
            \end{pmatrix}
            =\begin{pmatrix}
             \frac{\mathbf{v}_\bot-\mathbf{U}}{\varepsilon}\\\mathbf{v}_\|
            \end{pmatrix},
 \label{eq:changement_varialble}
\end{equation}
with
\begin{equation*}
 \mathbf{U}=c\left(\frac{\mathbf{E}\times\mathbf{B}}{B^2}\right)_\bot.
\end{equation*}
Next, we  decompose $\mathbf{E}$ into components along $\mathbf{B}$ and perpendicular to  $\mathbf{B}$ :
\begin{equation*}
 \mathbf{E}=\mathbf{E}_\bot+\mathbf{E}_{\|}e_z.
\end{equation*}
Substituting this expression in~\eqref{eq:changement_varialble}, it yields
\begin{equation*}
 \mathbf{E}+\frac{\mathbf{v}\times\mathbf{B}}{c}=\mathbf{E}_{\|}\mathbf{e}_z+\frac{\varepsilon}{c}\mathbf{w}\times\mathbf{B}.
\end{equation*}
Then we introduce a new distribution function $g$, such that
\begin{equation*}
g(t,\mathbf{x},\mathbf{w}) =  f(t,\mathbf{x},\mathbf{v})
\end{equation*}
for which we get that it is solution to
\begin{eqnarray*}
 \frac{\partial {g}}{\partial t}\,\,+\,\,(\mathbf{U}\,+\,\varepsilon{\mathbf{w}_\bot})\cdot\nabla_{\mathbf{x}_\bot}{g}\,\,+\;\,{\mathbf{w}}_\|\,\partial_z{g}\,\,+\;\,\frac{1}{\varepsilon}\,\mathbf{U}_{\bot}^T\,\nabla_{\mathbf{x}_\bot}\mathbf{U}\,\,\nabla_{{\mathbf{w}}_\bot}{g}\,+\,\frac{e}{m}E_{\|}\,\partial_{{\mathbf{w}}_\|}{g}
\\
 \,
 \,+\,\,{\mathbf{w}}_\bot^T\,\nabla_{\mathbf{x}_\bot}\,\mathbf{U}\,\nabla_{{\mathbf{w}}_\bot}{g}\,\,+\,\,\frac{e}{cm}\,{\mathbf{w}}\times\mathbf{B}\cdot\,\nabla_{{\mathbf{w}}}{g}\,\,\,=\,\,\,
 0.
\end{eqnarray*}
Finally, we   integrate the previous equation in velocity field ${\mathbf{w}}_\bot\in\mathbb{R}^2$
\begin{equation*}
     \frac{\partial \tilde{f}^\varepsilon}{\partial t}+\mathbf{U}\cdot\nabla_{\mathbf{x}_\bot}\tilde{f}^\varepsilon+{\mathbf{w}}_\|\partial_z\tilde{f}^\varepsilon+\frac{e}{m}E_{\|}\partial_{{\mathbf{w}}_\|}\tilde{f}^\varepsilon=-\varepsilon\int_{\mathbb{R}^2}{\mathbf{w}_\bot}\cdot\nabla_{\mathbf{x}_\bot}{g}\,d{\mathbf{w}}_\bot,
\end{equation*}
where  $\bar{f}^\varepsilon=\int_{\mathbb{R}^2}{g}\,d{\mathbf{w}}_\bot$.
By passing formally to the limit  $\varepsilon\to0$, we  obtain the  drift-kinetic model
  \begin{equation}
   \frac{\partial \bar{f}}{\partial t}\,+\,\mathbf{U}\cdot\nabla_{\mathbf{x}_\bot}\bar{f}\,+\,{\mathbf{w}}_\|\partial_z\bar{f}+\frac{e}{m}E_{\|}\partial_{{\mathbf{w}}_\|}\bar{f}\,=\,0.
   \label{eq:driftkinetic_bis}
  \end{equation}
  
  On the other hand, the self-consistent potential $\phi$  is solution
  to the quasi-neutrality equation~\cite{bibCS}
  \begin{equation}
   -\nabla_\bot\cdot\left(\frac{\rho _0}{B\omega_c}\nabla_\bot\phi\right)+\frac{e\, \rho_0}{T_e}(\phi-\bar{\phi})=\rho-\rho_0,
   \label{eq:quasi_neutre}
  \end{equation}
where $\omega_c=eB/m_i$ is the ion cyclotron frequency, and $T_e$ and
$\rho_0$ are the electron temperature and density profiles
respectively which only depend on $\mathbf{x}_\bot$. 
The ion density profile is given by 
$$
\rho(\mathbf{x},t)=\int_{\mathbb{R}}f(\mathbf{x},v_{\|},t)dv_{\|}
$$
and $\bar{\phi}$ represents the average on the magnetic field lines,
that is, 
$$
\bar{\phi}\,=\,\frac{1}{L_z}\,\int_0^{L_z}\phi dz,
$$ with $L_z$ the  length in the $z$ variable.

Dropping the $\bar{\,}$ over the various quantities and replacing $w$ by $v$, the non-dimensional form of  the  drift-kinetic model can be written as

\begin{equation}
\label{DK4D}
\left\{
\begin{array}{l}
 \displaystyle\frac{\partial {f}}{\partial
   t}+\mathbf{U}\cdot\nabla_{\mathbf{x}_\bot}{f}+{v}_\|\partial_z{f}+E_{\|}\partial_{{v}_\|}{f}=0,
\\[3mm]
 \displaystyle\mathbf{U}=\frac{\mathbf{E}\times\mathbf{B}}{B^2},\quad
 \mathbf{E}=-\nabla\phi,
\\[3mm]
\displaystyle -\nabla_\bot\cdot\left(\frac{\rho _0(\mathbf{x}_\bot)}{B}\nabla_\bot\phi\right)+\frac{\rho_0(\mathbf{x}_\bot)}{T_e(\mathbf{x}_\bot)}(\phi-\bar{\phi})=\rho-\rho_0.
\end{array}\right.
\end{equation}

The following proposition shows some properties of the drift-kinetic model when we ignore the effect of boundary conditions:
\begin{proposition}
\label{prop:DriftKinetic}
 Let us consider $\Omega=\RR^3$ or the torus. Then the  drift-kinetic model (\ref{DK4D}) verifies the following properties :
\begin{enumerate}
 \item  If $f$ is smooth, we have the maximum principle
 \begin{equation*}
   0\leq
   f(t,\mathbf{x},v_{\|})\leq\max_{\mathbf{x},v_{\|}}(f(0,\mathbf{x},v_{\|})),\quad
   t\,\geq\; 0.  
 \end{equation*}
 \item $L^p$ norm conservation, for $1\leq p\leq\infty$
 \begin{equation*}
  \frac{d}{dt}\int_{\RR}\int_{\Omega} |f(t,
  \mathbf{x},v_{\|})|^p\,d\mathbf{x}\,dv_\|\,\,=\,\,0,\quad t\geq 0.
   \end{equation*}
 \item Kinetic entropy conservation
 \begin{equation*}
  \frac{d}{dt} \int_{\RR}\int_{\Omega}
  f\ln|f|\,d\mathbf{x}\,dv_\| \,\,=\,\,0,\quad t\geq 0.
   \end{equation*}
\item Energy conservation
\begin{equation}
 \frac{d}{dt}\left(\int_{\RR}\int_{\Omega}  (f-f_M)v^2_\|   d\mathbf{x}dv_\|+\int_{\Omega}\phi(\rho-\rho_0)d\mathbf{x}\right)=0.
 \label{eq:energy_DK}
\end{equation}
\end{enumerate}
\end{proposition}

For practical applications,  this model has to be supplemented with
suitable boundary conditions when considering a domain of the form
\begin{equation*}
 \Omega=\{(x,y,z)\in\mathbb{R}^3 : (x,y)\in D, 0\leq z\leq L_z\},
\end{equation*}
with $D$ a two dimensional domain. We assume that the electric potential is vanishing at the boundary $\partial D$
\begin{equation}
 \phi(\mathbf{x})=0,\quad\mathbf{x}\in\partial D\times[0,L_z],
 \label{eq:DK_BC1}
\end{equation}
and the distribution function is given by 
\begin{equation}
 f(\mathbf{x}_\bot,z,v\|)=f_M(\mathbf{x}_\bot,z,v\|),\quad\mathbf{x}\in\partial
 D\times[0,L_z], \,{\rm with }\, \mathbf{U}(\mathbf{x})\cdot
 \mathbf{n}_\mathbf{x} \geq 0,
 \label{eq:DK_BC2}
\end{equation}
where $f_M$ is a stationary solution to the drift-kinetic model (\ref{DK4D}). 
Furthermore, periodic boundary conditions are assumed  for the potential and the distribution
function in the $z$ direction
\begin{equation*}
 \phi(\mathbf{x}_\bot,0)=\phi(\mathbf{x}_\bot,L_z), \quad
 f(\mathbf{x}_\bot,0,v_\|)=f(\mathbf{x}_\bot,L_z,v_\|),\quad\mathbf{x}_\bot\in
 D,\quad v_\|\in\RR.
\end{equation*}

Finally, we can derive the guiding centre model from (\ref{DK4D}) by
integrating the equation with
respect to $(z,v_\|)$. We get that the reduced density $\bar\rho :
\mathbf{x}_\bot \mapsto \RR$ is solution to the guiding centre system of
equations
\begin{equation}
\label{CG2D}
\left\{
\begin{array}{l}
 \displaystyle\frac{\partial {\bar\rho}}{\partial
   t}+\mathbf{U}\cdot\nabla_{\mathbf{x}_\bot}{\bar\rho} =0,
\\[3mm]
 \displaystyle\mathbf{U}=\frac{\mathbf{E}\times\mathbf{B}}{B^2}, \quad \mathbf{E}=-\nabla\phi,
\\[3mm]
\displaystyle -\nabla_\bot\cdot\left(\frac{\rho _0(\mathbf{x}_\bot)}{B}\nabla_\bot\phi\right)\,=\,\bar\rho-\rho_0.
\end{array}\right.
\end{equation}

In this case, the solution $\bar \rho$ satisfies the following
properties 

\begin{proposition}
\label{prop:CG2D}
 Let us consider the two dimensional domain $D=\RR^2$ or the torus.  Then the  guiding centre model (\ref{CG2D}) verifies the following properties :
\begin{enumerate}
 \item  If $\bar\rho$ is smooth, we have the maximum principle
 \begin{equation*}
   0\leq \bar\rho(t,\mathbf{x}_\bot)\leq\max_{\mathbf{x}_\bot}(\bar\rho(0,\mathbf{x}_\bot).  
 \end{equation*}
 \item $L^p$ norm conservation, for $1\leq p\leq\infty$
 \begin{equation*}
  \frac{d}{dt}\int_{D} |\bar\rho(t, \mathbf{x}_\bot)|^p
  \,d\mathbf{x}_\bot\,=\,0.
   \end{equation*}
\item Energy conservation
\begin{equation}
 \frac{d}{dt}\int_{D}  \bar\rho \,\phi\,d\mathbf{x}_\bot \,=\,0.
 \label{eq:energy_CG}
\end{equation}
\end{enumerate}
\end{proposition}

For practical applications, we assume that the electric potential is vanishing at the boundary $\partial D$
\begin{equation}
 \phi(t,\mathbf{x}_\bot)=0,\quad\mathbf{x}_\bot\in\partial D,\quad
 t\geq 0.
 \label{eq:CG_BC}
\end{equation}

\section{Numerical schemes}
\label{sec:schemes}
\setcounter{equation}{0}
In this section we present the hybrid method based on the Hermite WENO
reconstruction already proposed in \cite{YF15}. On the one hand, we
apply a semi-Lagrangian method for a general transport equation
written in a non conservative form. On the other hand, we apply a
finite difference method, which enforces the conservation of mass when
the equation is written in the conservative form.  These methods are
coupled with the inverse Lax-Wendroff procedure to discretize
accurately boundary conditions in an arbitrary $2D$ geometry.
Finally in the subsection \ref{sec:poisson},   we discretize the
Poisson equation for the electrical potential (\ref{CG2D}).


\subsection{Hermite WENO reconstruction for semi-Lagrangian methods}
\label{sec:HSL}
\setcounter{equation}{0}

We briefly remind the high order Hermite interpolation coupled with a
weight essentially non-oscillatory (HWENO) reconstruction for
semi-Lagrangian methods.  The semi-Lagrangian method   becomes a
classical method for the numerical solution of the Vlasov equation
because of its high accuracy and its small
dissipation~\cite{bibKnorr,bibSR}.  For a given $s\in\RR^+$, the differential system 
\begin{equation*}
 \left\{
 \begin{array}{l}
 \displaystyle \frac{d\mathbf{X}}{dt}= \mathbf{A}(t,\mathbf{X}),\\[3mm]
 \mathbf{X}(s)=\xx,
 \end{array}
 \right.
\end{equation*}
is associated to the transport equation
\begin{equation}
\label{eq:transport_vlasov}
\frac{\partial f}{\partial t} \,+\,
\mathbf{A}(t,\mathbf{X})\,\nabla_\mathbf{X} f = 0.
\end{equation}
 We denote its solution by $\mathbf{X}(t;s,\xx)$.
The backward semi-Lagrangian method is decomposed into two steps for computing the  function $f^{n+1}$ at time $t_{n+1}$ from the  function $f^{n}$ at time $t_{n}$ :
\begin{enumerate}
 \item For each mesh point $\xx_i$ of phase space, compute  the backward characteristic $\mathbf{X}(t_n;t_{n+1},\xx_i)$, the value of the characteristic at time $t_n$ who is equal to $\xx_i$ at time $t_{n+1}$.
 \item As the function $f$ of transport equation verifies
 \begin{equation*}
  f^{n+1}(\xx_i)=f^n(\mathbf{X}(t_n;t_{n+1},\xx_i)),
 \end{equation*}
we obtain the value of $f^{n+1}(\xx_i)$ by computing $f^n(\mathbf{X}(t_n;t_{n+1},\xx_i))$ by interpolation, since $\mathbf{X}(t_n;t_{n+1},\xx_i)$ is not usually a mesh point.
\end{enumerate}

We apply a third order Hermite interpolation coupled with a weighted
essentially non-oscillatory procedure,  such that  it is accurate for
smooth solutions and it removes  spurious oscillations around
discontinuities or high frequencies which cannot be solved on a fixed
mesh. Consider a uniform mesh $(x_i)_{i}$ of the computational domain
and assume that the values of the distribution function $(f_i)_i$ and
its derivative $(f'_i)_i$ are known at the grid points.  We define two quadratic polynomials in the interval $I_i$ :
$$
\left\{\begin{array}{l}
 \displaystyle h_l(x)\,=\,f_i+\frac{f_{i+1}-f_i}{\Delta x}(x-x_i)+\frac{(f_{i+1}-f_i)-\Delta x {f}'_i}{\Delta x^2}(x-x_i)(x-x_{i+1}),\label{eq:hweno:pl}\\[3mm]
 \displaystyle  h_r(x)\,=\,f_i+\frac{f_{i+1}-f_i}{\Delta x}(x-x_i)+\frac{\Delta x {f}'_{i+1}-(f_{i+1}-f_i)}{\Delta x^2}(x-x_i)(x-x_{i+1}).\label{eq:hweno:pr}
\end{array}\right.
$$
The  polynomial of degree 2 $h_l$ verifies 
\begin{equation*}
 h_l(x_i)=f_i,\quad h_l(x_{i+1})=f_{i+1},\quad h'_l(x_i)= f'_i,
\end{equation*}
while  $h_r$ verifies 
\begin{equation*}
 h_r(x_i)=f_i,\quad h_r(x_{i+1})=f_{i+1},\quad h'_r(x_{i+1})= f'_{i+1}.
\end{equation*}

The idea of WENO reconstruction is now to use the cubic polynomial when function $f$ is smooth, otherwise, we use the less oscillatory  polynomial of degree 2 between $h_l$ or $h_r$. Thus, let us define  $H_3$ as follows
\begin{equation*}
 H_3(x)\,\,=\,\,w_l(x)\,h_l(x)\,\,+\,\,w_r(x)\,h_r(x),
\end{equation*}
where $w_l$ and $w_r$ are WENO weights. To determine these WENO weights, we follow the strategy given in~\cite{bibJS} and first define smoothness indicators by integration of the first and second derivatives of  $h_{l}$ and $h_{r}$ on the interval $I_i$ :
$$
\left\{\begin{array}{l}
 \displaystyle \beta_l \,=\,\int_{x_i}^{x_{i+1}}\Delta x(h_{l}')^2+\Delta x^3(h_{l}'')^2dx\,=\,(f_i-f_{i+1})^2+\frac{13}{3}((f_{i+1}-f_i)-\Delta x {f}'_i)^2,
\\
\,
\\
 \displaystyle \beta_r\,=\;\int_{x_i}^{x_{i+1}}\Delta x(h_{r}')^2+\Delta x^3(h_{r}'')^2dx\,=\;(f_i-f_{i+1})^2+\frac{13}{3}(\Delta x {f}'_{i+1}-(f_{i+1}-f_i))^2.
\end{array}\right.
$$
Then we set $w_l$ and $w_r$ as
\begin{equation*}
 w_l(x)=\frac{\alpha_l(x)}{\alpha_l(x)+\alpha_r(x)}\quad{\rm and}\quad w_r(x)=1-w_l(x),
 \end{equation*}
where
\begin{equation*}
 \alpha_l(x)=\frac{c_l(x)}{(\varepsilon+\beta_l)^2}\quad{\rm and}\quad \alpha_r(x)=\frac{c_r(x)}{(\varepsilon+\beta_r)^2}.
\end{equation*}
where $c_l=(x_{i+1}-x)/\Delta x$,  $c_r=1-c_l$ and $\varepsilon=10^{-6}$ to avoid the denominator to be zero.

Observe that when the function  $f$ is smooth, the difference between
$\beta_l$ and $\beta_r$ becomes small  and the weights  $w_l(x)\approx
c_l(x)$ and $w_r(x)\approx c_r(x)$.  Otherwise, when the smoothness
indicator $\beta_s$, $s=l,r$ blows-up, then the parameter $\alpha_s$
and the weight $w_s$ goes to zero, which
yields
\begin{equation}
  w_l(x)\approx 1, w_r(x)\approx 0\quad \text{or}\quad w_l(x)\approx 0, w_r(x)\approx 1.
  \label{eq:weno_approx2}
 \end{equation}
  Finally, let us mention that the following fourth-order centred finite difference formula is used to approximate the first derivative at the grid point $x_i$
 \begin{equation}
  f'_i=\frac{1}{12\Delta x}(f_{i-2}-8f_{i-1}+8f_{i+1}-f_{i+2}).
  \label{eq:4order_dxf}
 \end{equation}

\subsection{Hermite WENO reconstruction for conservative finite difference methods}
\label{sec:HFD}
\setcounter{equation}{0}
When the velocity $\mathbf{A}$ is not constant in (\ref{eq:transport_vlasov}), the semi-Lagrangian method is not conservative even when $\Div \mathbf{A}=0$, hence mass is no longer conserved and the long time behavior of the numerical solution may be wrong even for small time steps. Therefore, high order conservative methods may be more appropriate even if they are restricted by a CFL condition.

In this section, we extend Hermite WENO reconstruction for computing numerical flux of finite difference method.
Suppose that $\{f_i\}_{1\leq i\leq N}$  is approximation of $f(x_i)$. 
We look for $\{\hat f_{i+1/2}\}_{0\leq i\leq N}$ such that the flux difference approximates the derivative $f'(x)$ to $k$-th order accuracy :
\begin{equation*}
 \frac{ \hat f_{i+1/2} - \hat f_{i-1/2} }{\Delta x} = f'(x) + \mathcal{O} ( \Delta x^k ).
\end{equation*}
To approximate the flux $\hat f_{i+1/2} $, we define  a piecewise
polynomial $G$ such that it is exactly known on a set of points
$x=x_{k+1/2}$, $k=i-l,\ldots,i+r$.
\begin{equation*}
G(x_{k+1/2} ) \,=\,  G_{k+1/2}  \, =\, \Delta x \sum_{j=-\infty} ^k
f_j, \quad i-l\leq k\leq i+r. 
\end{equation*}
Thus, given the point values $\{f_i\}$,  we can compute  $G(x)$ by an interpolation method and therefore deduce the
numerical flux by 
\begin{equation}
 \hat f_{i+1/2}   = \left.\frac{d G}{d x}\right|_{x=x_{i+1/2}}.
\end{equation}

Now to interpolate the  function $G(x)$, we  apply a high order Hermite WENO scheme and outline the procedure of reconstruction only for the fifth order accuracy case.

The aim is to construct an approximation of the flux  $f^{-}_{i+1/2}$ by the Hermite polynomial of degree five together with a WENO reconstruction from point values $\{f_i\}$ :
\begin{enumerate}
 \item We construct the Hermite polynomial $H_5$ such that
$$
H_5(x_{i+j+1/2}) = G_{i+j+1/2},\, j=-2,\,-1,\,0,\,1,\quad H'_5(x_{i+j+1/2}) =  G'_{i+j+1/2},\, j=-2,\,1,
$$
\item We construct  cubic reconstruction polynomials $H_l(x)$, $H_c(x)$, $H_r(x)$ such that :
 $$
\left\{
\begin{array}{ll}
 H_l(x_{i+j+1/2}) = G_{i+j+1/2},\, j=-2,\,-1,\,0, & H'_l(x_{i-3/2}) =  G'_{i-3/2},\\[4mm]
 H_c(x_{i+j+1/2}) = G_{i+j+1/2},\, j=-2,\,-1,\,0,\,1,& \, \\[4mm]
 H_r(x_{i+j+1/2}) = G_{i+j+1/2},\, j=-1,\,0,\,1, & H'_r(x_{i+3/2}) =  G'_{i+3/2},\\[4mm]
  \end{array}\right.
$$
where $G'_{i+1/2}$ is the sixth order centered approximation of first derivative
\begin{equation*}
  G'_{i+1/2} = \frac{1}{60}[(u_{i+3} + u_{i-2})   -  8(u_{i+2} + u_{i-1})   +  37(u_{i+1} + u_{i})].
\end{equation*}
 Let us denote by $h_l(x)$, $h_c(x)$, $h_r(x)$, $h_5(x)$  the first derivatives of $H_l(x)$, $H_c(x)$, $H_r(x)$, $H_5(x)$  respectively. By evaluating  $h_l(x)$, $h_c(x)$, $h_r(x)$, $h_5(x)$ at $x=x_{i+1/2}$, we obtain
$$
 h_5(x_{i+1/2}) \,=\, \frac{ -8f_{i-1} + 19f_{i} + 19f_{i+1}   + 3G'_{i-3/2} - 6G'_{i+3/2}}{27}
$$
and
$$\left\{\begin{array}{l}
 \displaystyle  h_l(x_{i+1/2}) \,=\  -2\,f_{i-1} \,+\, 2f_{i} \,+\, G'_{i-3/2} , \\[3mm]
 \displaystyle   h_c(x_{i+1/2}) \,=\, \frac{-f_{i-1} \,+\, 5\,f_i + 2\,f_{i+1}}{6},\\[3mm]
 \displaystyle   h_r(x_{i+1/2}) \,=\, \frac{ f_{i} \,+\, 5\,f_{i+1}  \,-\, 2\,G'_{i+3/2}}{4}.
\end{array}\right.
$$

\item We evaluate the smoothness indicators $\beta_l$, $\beta_c$,
  $\beta_r$, which measure the smoothness of $h_l(x)$, $h_c(x)$,
  $h_r(x)$ on the cell $[x_i, x_{i+1}]$
\begin{eqnarray*}
 \beta_l &=&\int_{x_i}^{x_{i+1}}  \Delta x (h'_l(x))^2 +  \Delta x^3 (h''_l(x))^2 dx\\[3mm]
 &=& l_1^2 \,\,+\,\, 3\, l_1\,l_2 \,\,+\,\, \frac{75}{16} \,l_2^2, \quad \textrm{with } l_1 \,=\, u_i - u_{i-1}, \quad l_2 \,=\, -3u_{i-1} + u_i + 2G'_{i-3/2},     \\[3mm]
 \beta_c &=&\int_{x_i}^{x_{i+1}}  \Delta x (h'_c(x))^2 +  \Delta x^3 (h''_c(x))^2 dx\\[3mm]
 &=& c_1^2 \,\,+\,\, 2\, c_1\,c_2 \,\,+\,\, \frac{25}{12} \,c_2^2,
 \quad \textrm{with } c_1 \,=\, u_i - u_{i-1}, \quad c_2 \,=\, u_{i-1} - 2u_i + u_{i+1},\\[3mm]
 \beta_r &=&\int_{x_i}^{x_{i+1}}  \Delta x (h'_r(x))^2 +  \Delta x^3 (h''_r(x))^2 dx\\[3mm]
 &=& r_1^2 \,\,+\,\, \frac{39}{16} \,r_2^2,
 \quad \textrm{with } r_1 \,=\, u_{i+1} - u_{i}, \quad r_2 \,=\, u_{i} - 3 u_{i+1} + 2G'_{i+3/2}.
\end{eqnarray*}

\item We compute the nonlinear weights based on the smoothness indicators
$$
\left\{
\begin{array}{ll}
 \displaystyle w_l = \frac{\alpha_l}{\alpha_l + \alpha_c + \alpha_r }, &  \displaystyle\alpha_l = \frac{c_l}{(\varepsilon + \beta_l)^2}, \\[5mm]
  \displaystyle w_c = \frac{\alpha_c}{\alpha_l + \alpha_c + \alpha_r }, &  \displaystyle\alpha_c = \frac{c_c}{(\varepsilon + \beta_c)^2}, \\[5mm]
 \displaystyle  w_r = \frac{\alpha_r}{\alpha_l + \alpha_c + \alpha_r }, &  \displaystyle\alpha_r = \frac{c_r}{(\varepsilon + \beta_r)^2}, 
\end{array}\right.
$$
where the coefficients $c_l=1/9$, $c_c=4/9$, $c_r=4/9$ are chosen to get fifth order accuracy for smooth solutions and the parameter $\varepsilon=10^{-6}$ avoids the blow-up of $\alpha_k$, $k=\{l,c,r\}$.

\item The flux $f^-_{i+1/2} $ is then computed as
\begin{equation*}
 f^-_{i+1/2} \,\,=\,\; w_l \,h_l(x_{i+1/2})  \,\,+\,\, w_c\, h_c(x_{i+1/2})  \,\,+\,\, w_r \,h_r(x_{i+1/2}). 
\end{equation*}

\end{enumerate}
The reconstruction to $f^+_{i+1/2}$ is mirror symmetric with respect
to $x_{i+1/2}$ of the above procedure.


\subsection{\bf Discretization of the Poisson  equation~\eqref{CG2D}-\eqref{eq:CG_BC}}

\label{sec:poisson}

We use  a classical five points finite difference approximation to
discretize  the Poisson equation~\eqref{CG2D}-\eqref{eq:CG_BC}.  So it remains to treat
the Dirichlet boundary conditions on $\partial D$.

\begin{figure}[h]
  \begin{center} 
    \includegraphics[width=10cm]{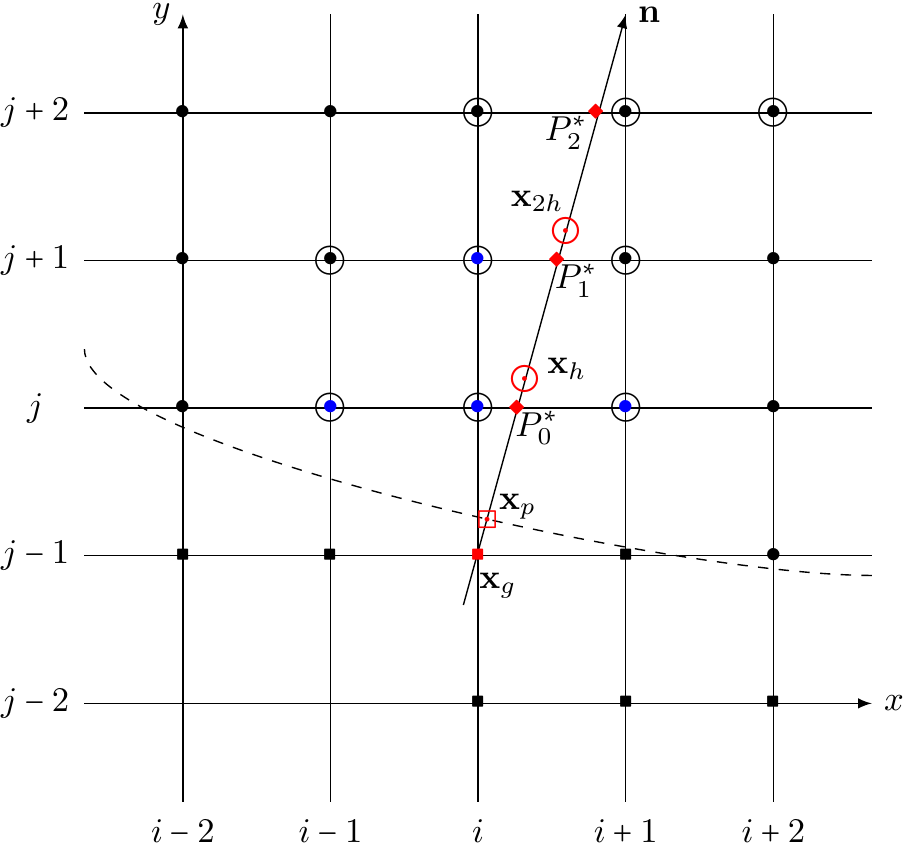}
\caption{\label{fig:2Ddomain}Spatially two-dimensional Cartesian mesh. $\bullet$ is interior point, $\filledsquare$ is ghost point, $\boxdot$ is the point at the boundary, $\largecircle$ is the point for extrapolation, the dashed line is the boundary.}
  \end{center}
\end{figure}

To discretize the Laplacian operator $\Delta_{\mathbf{x}_\bot}\phi$ near the physical boundary, some points of the usual five points finite difference formula can be located outside of interior domain. 
For instance,   Figure~\ref{fig:2Ddomain} illustrates the discretization stencil for  $\Delta_{\mathbf{x}_\bot}\phi$ at the point $(x_i,y_j)$.
We notice that the point $\mathbf{x}_g=(x_i,y_{j-1})$ is located outside of interior domain. Let us denote the approximation of $\phi$ at  the point $\mathbf{x}_g$ by $\phi_{i,j-1}$.
Thus $\phi_{i,j-1}$ should be extrapolated from the interior domain.

We extrapolate $\phi_{i,j-1}$ on the normal direction $\mathbf{n}$ 
\begin{equation}
 \phi_{i,j-1}=\bar{w}_p {\phi}(\mathbf{x}_p) + \bar{w}_h {\phi}(\mathbf{x}_h) + \bar{w}_{2h} {\phi}(\mathbf{x}_{2h}) ,
 \label{eq:normal_extrapolation}
\end{equation}
where $\mathbf{x}_p$ is the cross point of the normal $\mathbf{n}$ and the physical boundary $D$.
The points $\mathbf{x}_h$ and $\mathbf{x}_{2h}$ are  equal spacing on the normal $\mathbf{n}$, {\it i.e.} $h=|\mathbf{x}_p-\mathbf{x}_h|=|\mathbf{x}_h-\mathbf{x}_{2h}|$,
with $h=\min(\Delta x,\Delta y)$, $\Delta x$, $\Delta y$ are the space steps in the directions $x$ and $y$ respectively.
Moreover, $\bar{w}_p$, $\bar{w}_h$, $\bar{w}_{2h}$ are the extrapolation weights depending on the position of $\mathbf{x}_g$, $\mathbf{x}_p$, $\mathbf{x}_h$ and $\mathbf{x}_{2h}$.
In~\eqref{eq:normal_extrapolation}, ${\phi}(\mathbf{x}_p)$ is given by the  boundary condition~\eqref{eq:CG_BC},
whereas ${\phi}(\mathbf{x}_h)$, ${\phi}(\mathbf{x}_{2h})$ should be determined by interpolation.

For this, we first construct an interpolation stencil $\mathcal{E}$, composed of grid points of $D$. 
For instance, in Figure~\ref{fig:2Ddomain},  the inward normal $\mathbf{n}$ intersects  the grid lines $y=y_{j}$, $y_{j+1}$, $y_{j+2}$ at points $P^*_0$, $P^*_1$, $P^*_2$. 
Then we choose the three nearest points  of the cross point $P^*_l,\,l=0,1,2$, in each line, {\it i.e.} marked by a large circle. 
From these nine points, we can build a  Lagrange polynomial $q_2(\mathbf{x})\in\mathbb{Q}_2(\mathbb{R}^2)$.
Therefore,  we evaluate the polynomial $q_2(\mathbf{x})$ at $\mathbf{x}_h$ and $\mathbf{x}_{2h}$, {\it i.e.}
\begin{eqnarray*}
 {\phi}(\mathbf{x}_h)&=&\sum_{\ell=0}^8w_{h,\ell}{\phi}(\mathbf{x}_{\ell}),\\[3mm]
 {\phi}(\mathbf{x}_{2h})&=&\sum_{\ell=0}^8w_{2h,\ell}{\phi}(\mathbf{x}_{\ell}),
\end{eqnarray*}
with $\mathbf{x}_{\ell}\in \mathcal{E}$. 
We thus have that $\phi_{i,j-1}$ is approximated from the interior domain.

However, in some cases, we can not find a stencil of nine interior points. 
For instance, when the interior domain has small acute angle sharp, the normal $\mathbf{n}$ can not have three cross points $P^*_l,\,l=0,1,2$ in interior domain, or we can not have three nearest points  of the cross point $P^*_l,\,l=0,1,2$, in each line. 
In this case, we alternatively use a first degree polynomial $q_1(\mathbf{x})$ with a four points  stencil or even a zero degree polynomial $q_0(\mathbf{x})$ with  an one point  stencil. 
We can similarly construct the four points stencil or the one point stencil as above.


\section{Numerical simulations}
\label{sec:Num}
\setcounter{equation}{0}
In this section, we present numerical simulations. We first consider the 2D guiding-center model in a D-shaped domain. The steady state solution and perturbed solution are studied.
Then we perform the ion turbulence instability simulation with the 4D Drift-Kinetic model in a cylinder domain.
\subsection{\bf Test 1 : Steady state solution for the guiding center model in a D-shaped domain}
\begin{figure}
\begin{center}
 \begin{tabular}{ccc}
  \includegraphics[width=4.5cm]{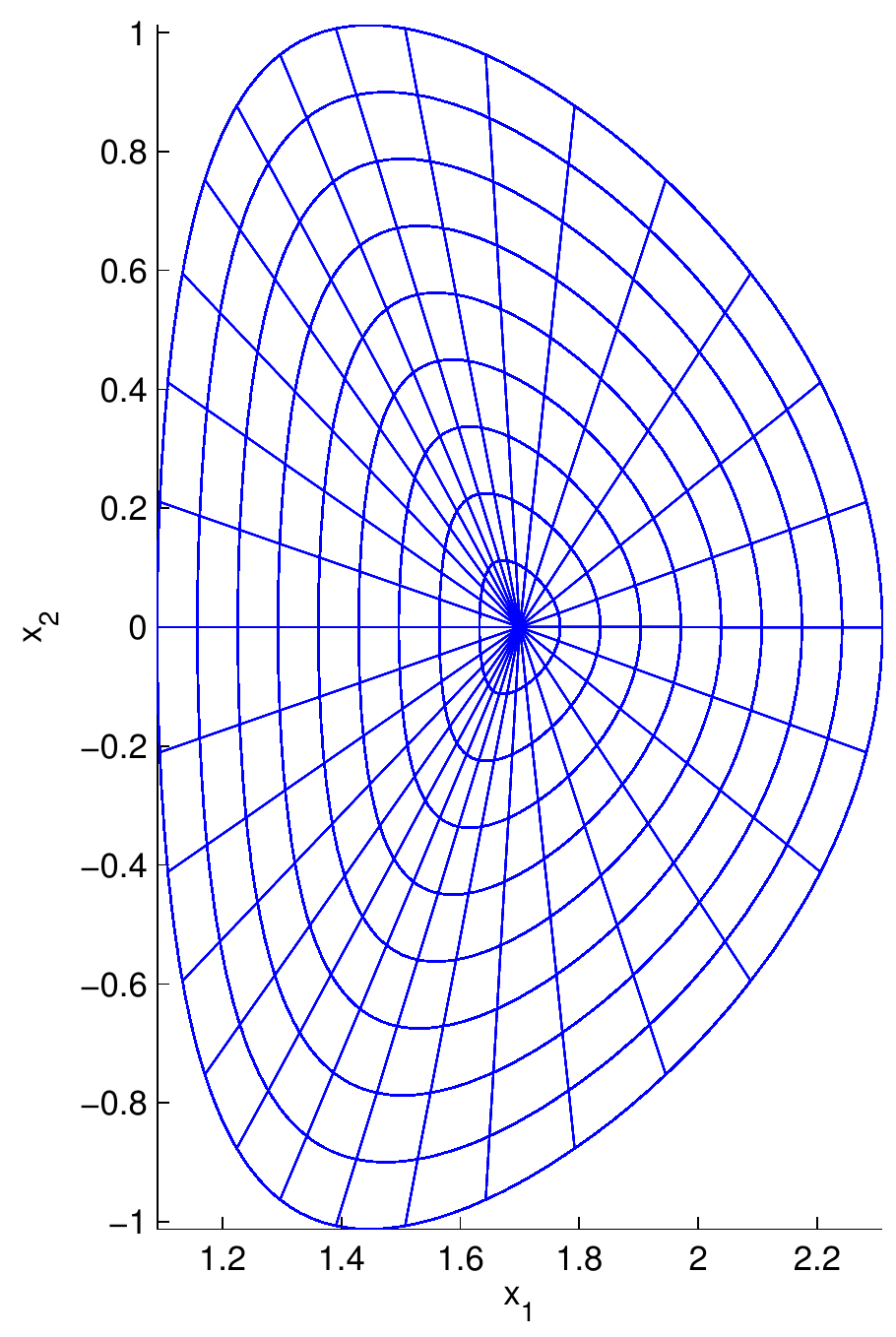} & $\quad\quad\quad$&   \includegraphics[width=4.5cm]{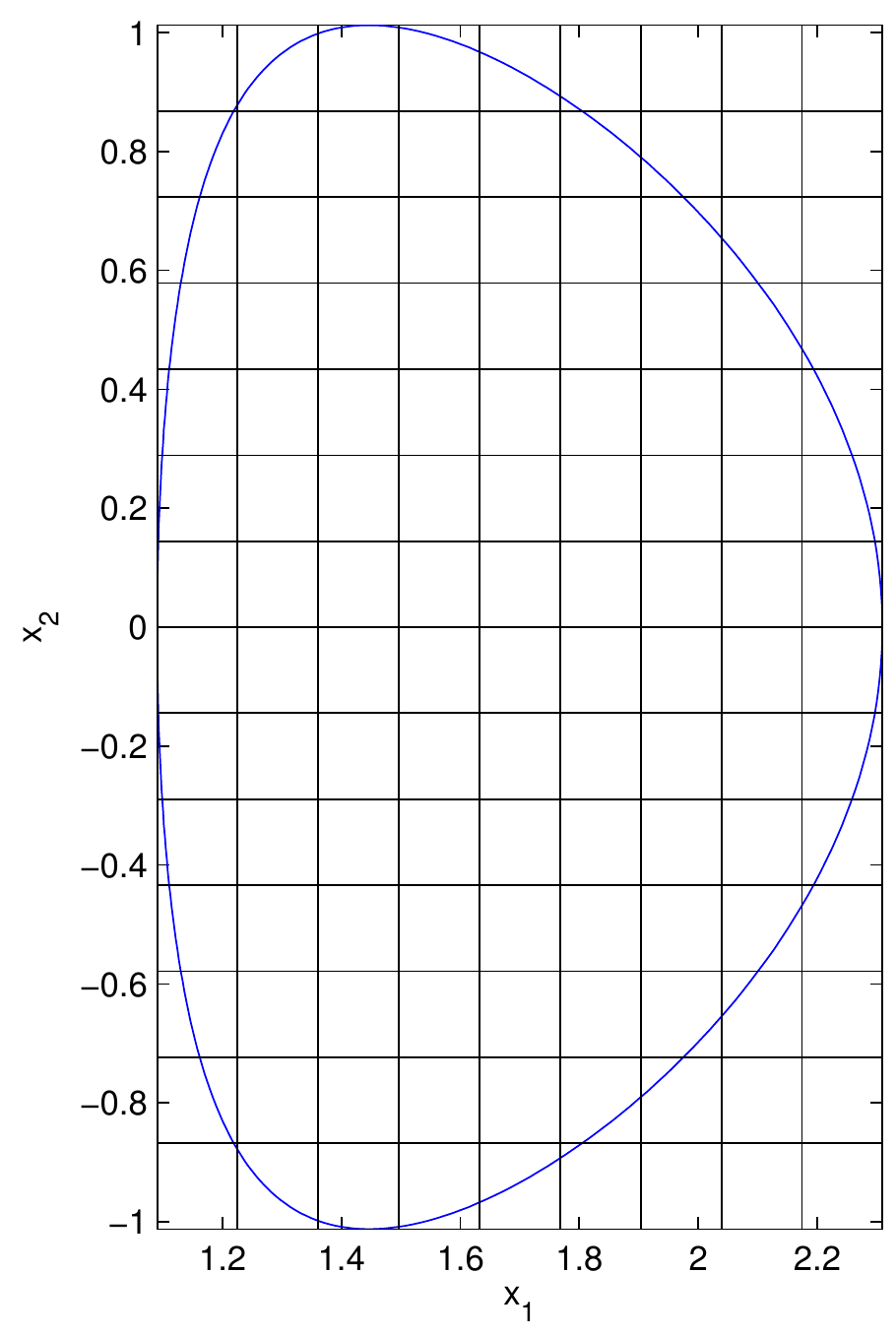} \\
  (a)                                       &  &  (b)
 \end{tabular}
\caption{\label{fig:Dshape}D-shaped domain. (a) Constant lines in coordinates $\xi = (\xi_1,\xi_2)$; (b) D-shaped domain embedded in Cartesian mesh.}
 \end{center}
\end{figure}

We consider now the 2D guiding-center model in a D-shaped domain $\Omega$ presented in Section IV of~\cite{bibmiller} and depicted in Figure~\ref{fig:Dshape}(a). The mapping $\mathbf{X}$ from curvilinear coordinates  $\mathbf{\xi} = (\xi_1,\xi_2)$ to physical coordinates $\mathbf{x} = ({x}_1, {x}_2)$ is given by

\begin{equation*}
  \begin{array}{l}
    x \, = \, 1.7 + [ 0.074(2\xi_1 - 1) + 0.536]\, \cos[2\pi\xi_2 + \arcsin(0.416)\sin(2\pi\xi_2)],\\[3mm]
    y\, = \, 1.66[ 0.074(2\xi_1 - 1) + 0.536]\sin(2\pi\xi_2),
  \end{array}
\end{equation*} 
for $-231/74 \, \leq \, \xi_1 \, \leq \, 1$, $0 \, \leq \, \xi_2 \, \leq \, 1$.

We now search a steady state solution for the guiding-center model in the D-shaped domain $\Omega$. We first notice that $\bar{\rho}(\phi)$, a function of $\phi$, is a solution of the guiding-center equation~\eqref{CG2D}. Then it remains to solve a nonlinear elliptic equation
\begin{equation}
  \left\{
  \begin{array}{ll}
    -\nabla_\bot \cdot \left( \frac{\rho_0}{B} \nabla_\bot \phi \right) = \bar{\rho}(\phi)- \rho_0 & \text{in } \Omega, \\[3mm]
    \phi = 0 & \text{on } \partial\Omega.
  \end{array}
\right.
\label{eq:nonlin_ellip}
\end{equation}
For a suitable function $\bar{\rho}$, we have a unique solution of equation~\eqref{eq:nonlin_ellip} :
\begin{proposition}
  Consider the equation~\eqref{eq:nonlin_ellip} with $\rho_0 = 1$, $B = 1$, $\bar{\rho}(\phi) = e^{-\phi}-1$. The function $\bar\rho$ is bounded on
  \begin{equation*}
    C \, = \, \{ \phi \, \in \, H^1_0(\Omega) : \phi \geq -\ln2 \} \subset H^1_0(\Omega).
  \end{equation*}
Then~\eqref{eq:nonlin_ellip} has a unique weak solution in $C$.
\end{proposition}
In the previous proposition, the existence of a positive solution is obtained by Schauder's fixed point theorem~\cite{bibEvans}, while the uniqueness is straight.

\begin{figure}
\begin{center}
 \begin{tabular}{cc}
  \includegraphics[width=7.5cm]{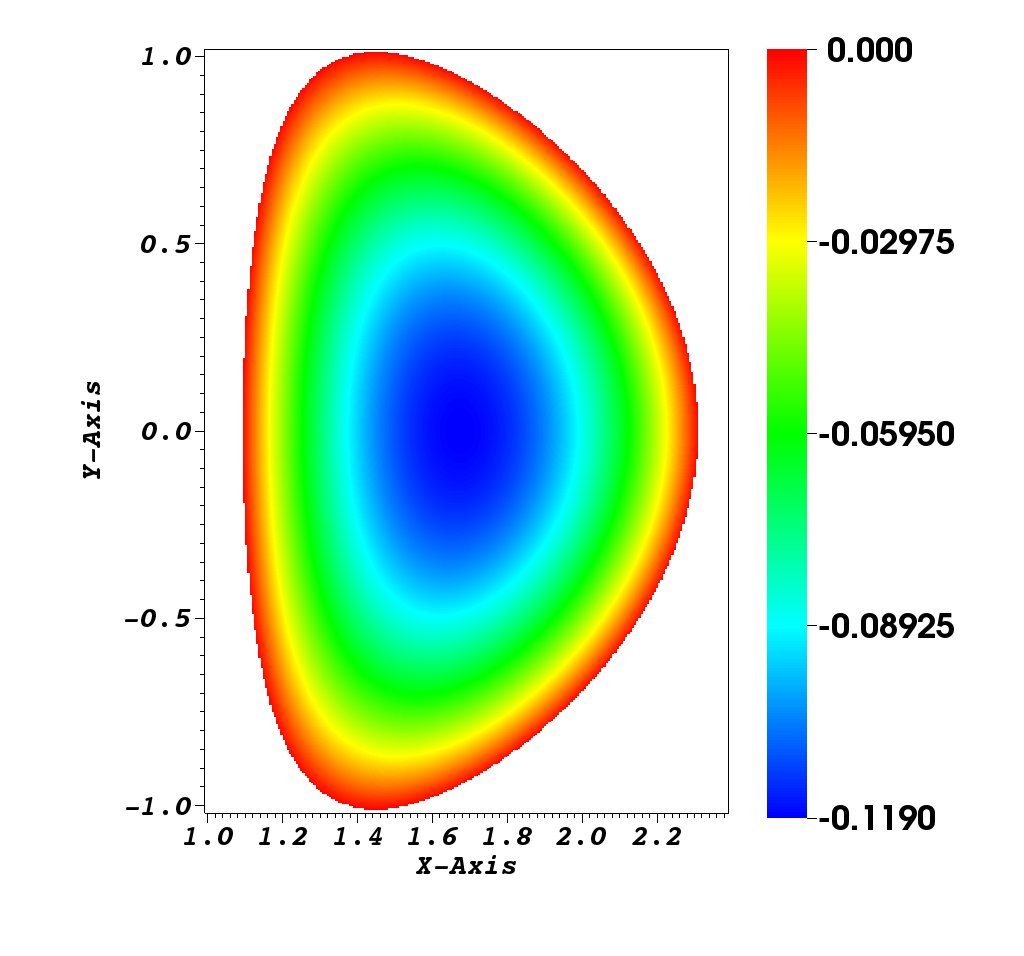} &    \includegraphics[width=7.5cm]{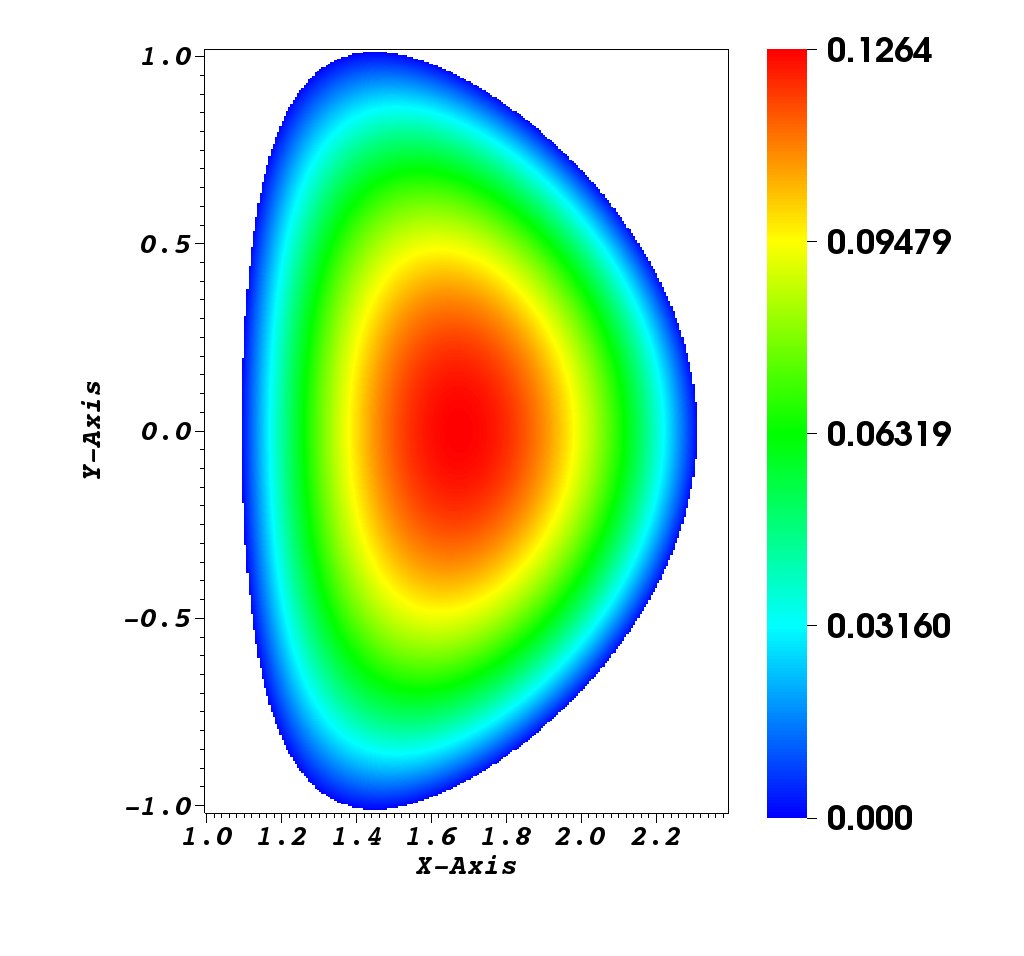} \\
  (a)   Potential $\phi_0$                                   &    (b) Density $\rho_0$
 \end{tabular}
\caption{\label{fig:Dshape_steady}A steady state solution of the guiding-center model~\eqref{CG2D} in D-shaped domain. Mesh size is $n_{x}\times n_{y} = 240\times 440$.}
 \end{center}
\end{figure}

\begin{figure}
\begin{center}
  \includegraphics[width=4.5cm]{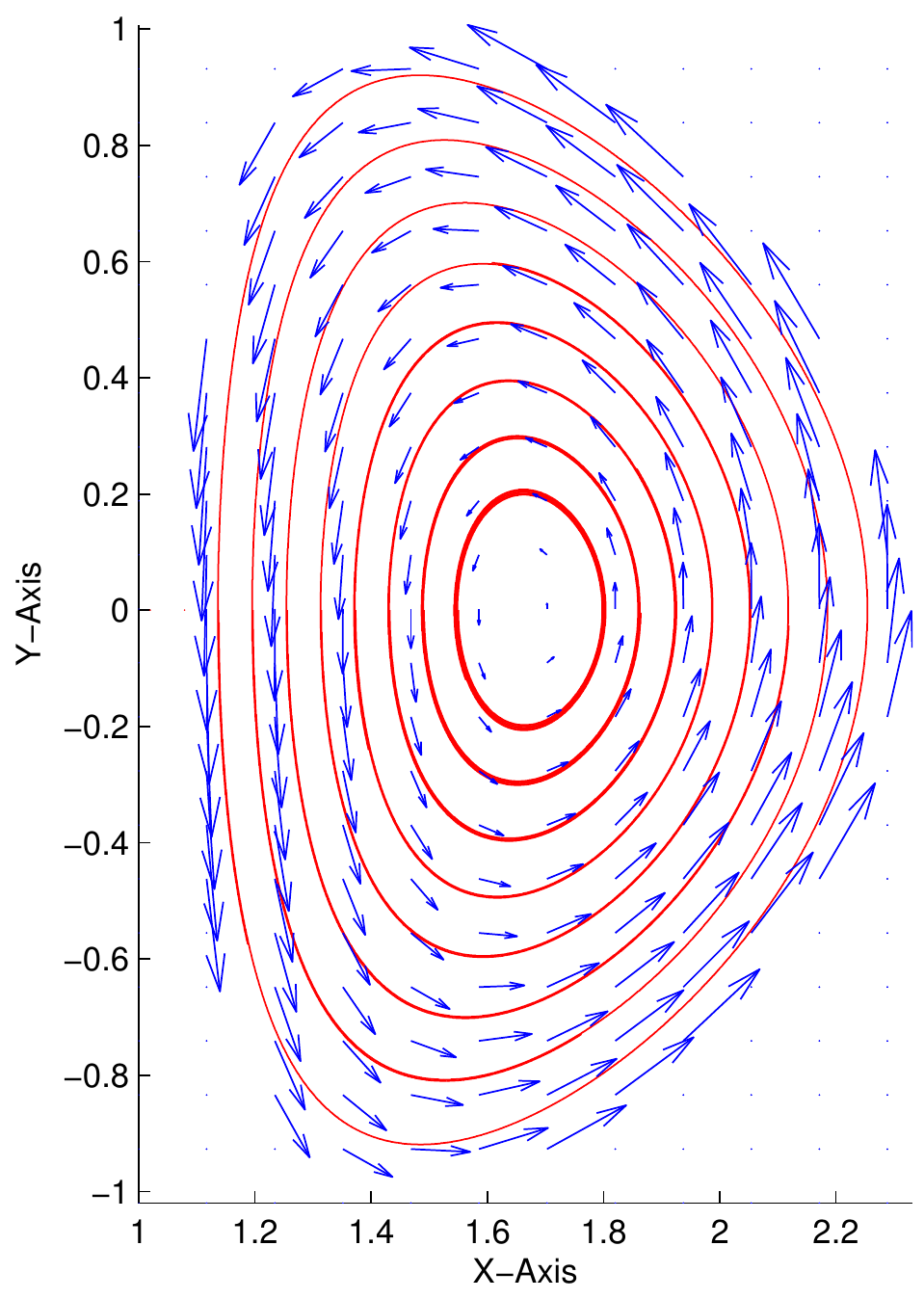}
\caption{\label{fig:Dshape_streamline}Streamline of the velocity field corresponding to steady state solution of guiding-center model~\eqref{CG2D} in D-shaped domain. }
 \end{center}
\end{figure}

We now solve~\eqref{eq:nonlin_ellip} in $\Omega$. Figure~\ref{fig:Dshape}(b) illustrates that the boundary $\partial\Omega$ is embedded in Cartesian mesh. Thus the numerical scheme presented in section~\ref{sec:poisson} should be applied, which gives us a nonlinear system for $\phi$. Then by applying a Newton method to this nonlinear system, we obtain a steady state solution $\phi_0$ of~\eqref{eq:nonlin_ellip} shown in Figure~\ref{fig:Dshape_steady}(a). Then by injecting $\phi_0$ into the guiding-center equation~\eqref{CG2D}, we get the steady state density $\bar\rho_0$ (see Figure~\ref{fig:Dshape_steady}(b)). Moreover, we plot the velocity field of steady state solution in Figure~\ref{fig:Dshape_streamline}. By comparing the streamline in Figure~\ref{fig:Dshape_streamline} and the constant line of  coordinates  $\mathbf{\xi} = (\xi_1,\xi_2)$, it is interesting to notice that these lines don't coincide, especially in the core of the D-shaped domain. The velocity is fast near the edge of the D-shaped domain but much slow in the core.

\begin{figure}
\begin{center}
 \begin{tabular}{cc}
  \includegraphics[width=7.5cm]{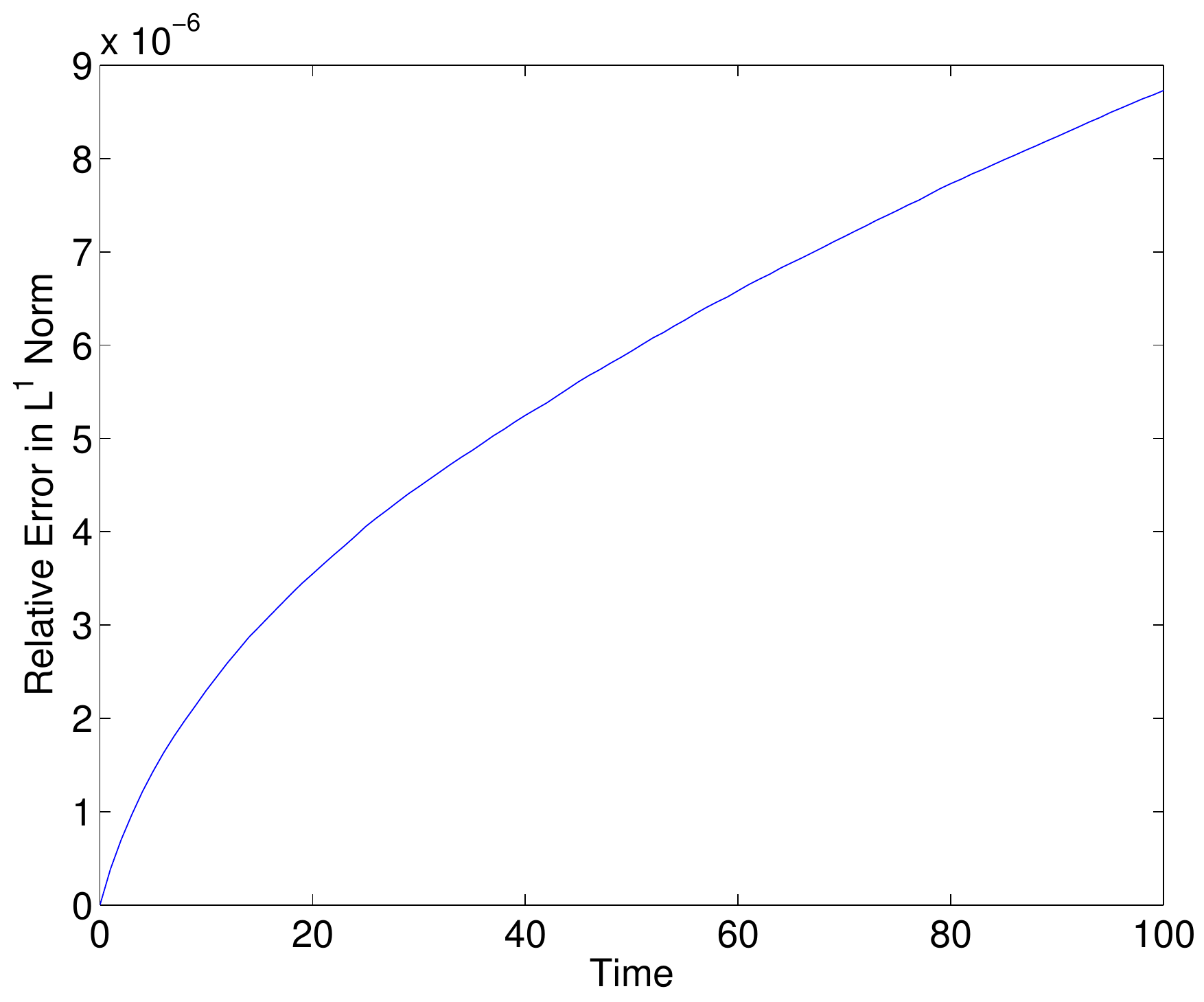} &    \includegraphics[width=7.5cm]{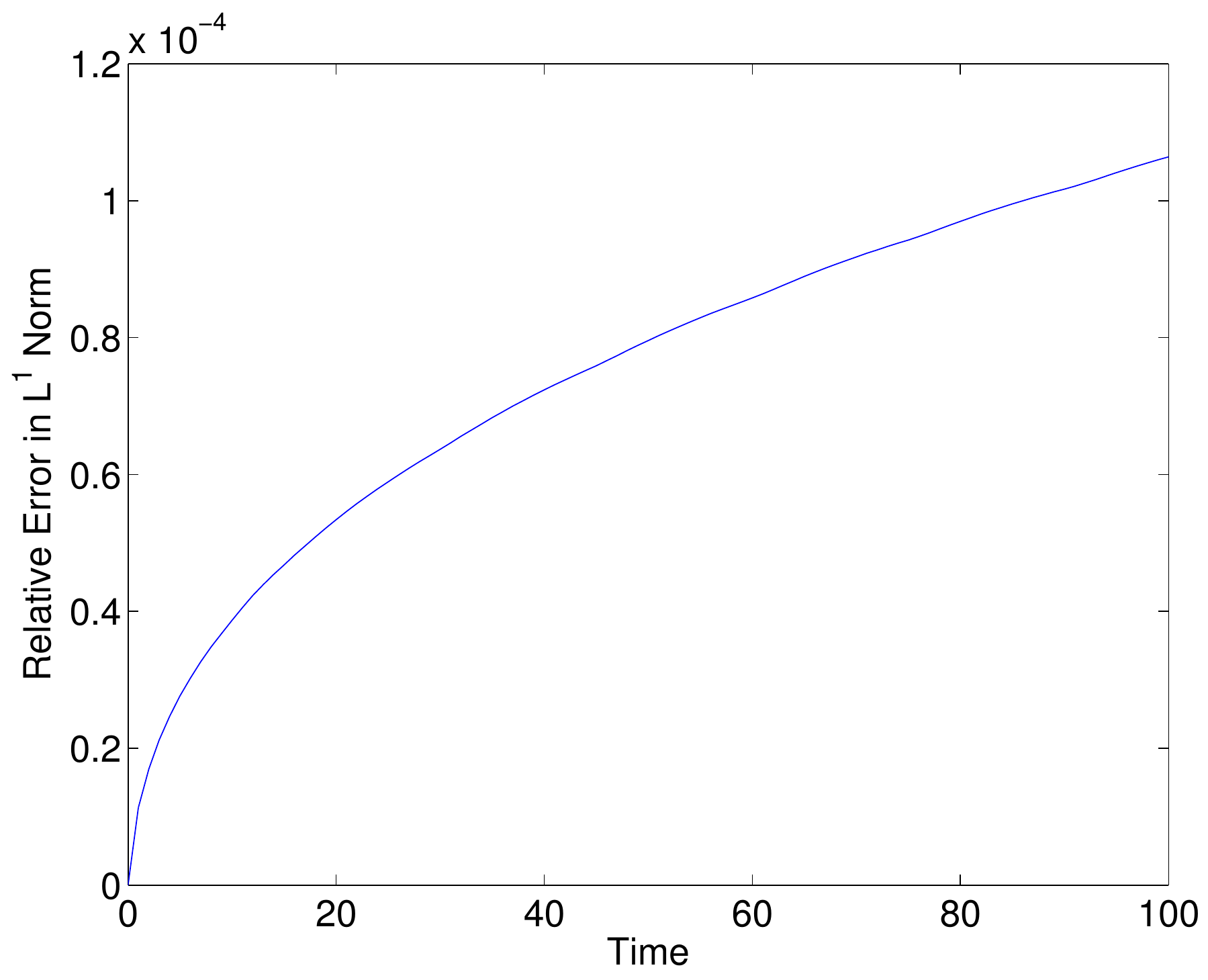} \\
  (a)   Potential                                    &    (b) Density
 \end{tabular}
\caption{\label{fig:Dshape_diff}Relative error in $L^1$ norm for potential $\phi_h$ and density $\rho_h$.}
 \end{center}
\end{figure}

Now we verify that $(\phi_0,\bar\rho_0)$ is the steady state solution by long time simulation. That is to take the pair $(\phi_0,\bar\rho_0)$ as an initial solution for  the guiding-center model~\eqref{CG2D}, then to compare the difference between $(\phi(t),\bar\rho(t))$ and  $(\phi_0,\bar\rho_0)$. We will measure these differences by a relative error as 
\begin{equation*}
  E(u(t)) \, = \, \frac{\|u(t) - u_0\|_1}{\|u_0\|_1},\, u = \phi,\,\bar\rho.
\end{equation*}
  The finite difference method with HWENO reconstruction is used for solving guiding-center equation~\eqref{CG2D}. The time step is taken to be $\Delta t = 0.001$. Figure~\ref{fig:Dshape_diff} presents the relative errors of the potential $\phi$ and the density $\bar\rho$. We observe that the solution remains steady for long time simulation with a relative error of magnitude of $10^{-4}$.

\subsection{\bf Test 2 : Numerical simulation of the guiding center model in a
  D shape}
Now we still consider the previous initial data $(\phi_0, \bar\rho_0)$ which is a stationary
solution of the guiding-center model, but perturb it of magnitude of  $\varepsilon$. 

In Figure~\ref{fig:Dshape_streamline}, we have noticed that the streamline is different with respect to the constant line of coordinates $(\xi_1, \xi_2)$. On the other hand, we denote
\begin{equation*}
  \mathbf{U} = \begin{pmatrix} U_{x} \\[3mm]
                               U_{y} \end{pmatrix}
            =  \begin{pmatrix}- \frac{\partial \phi}{\partial y} \\[3mm]
                              \frac{\partial \phi}{\partial x} \end{pmatrix}.
\end{equation*}
Then by the definition of streamline, we have
\begin{equation*}
  U_{y} - U_{x} \frac{d y}{d x } = 0,
\end{equation*}
which implies
\begin{equation*}
  \frac{d}{d x}\phi(x, y(x)) = \frac{\partial \phi}{\partial x} + \frac{\partial \phi}{\partial x} \frac{d y}{d x}
= U_{y} - U_{x} \frac{d y}{d x} = 0.
\end{equation*}
Thus, 
\begin{equation*}
  \phi(x, y) = const
\end{equation*}
represents a streamline, {\it i.e.} the steady state function $\bar\rho_0$ revolves along the isoline of potential $\phi_0$. 
In this test case, we perturb the function  $\bar\rho_0$ along the streamline, that is
\begin{equation*}
\bar\rho = \bar\rho_0 (1+\varepsilon\,\cos(2\pi k\xi_2)\,\exp(-2|\phi_0 - \phi_p|^2 / \varepsilon^4), 
\end{equation*}
with $\phi_p=-0.1$, $k=5$ and $\varepsilon=0.1$.

\begin{figure}
\begin{center}
 \begin{tabular}{cc}
  \includegraphics[width=7.5cm]{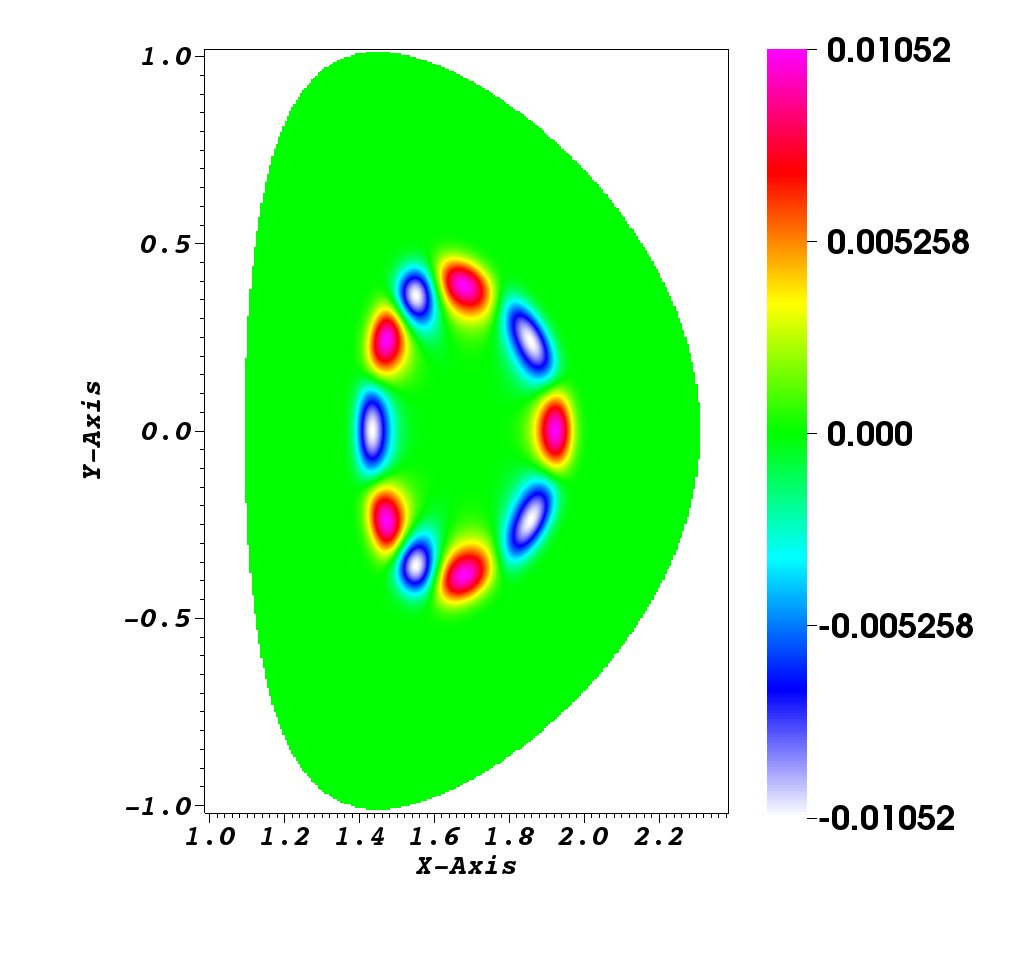} &    \includegraphics[width=7.5cm]{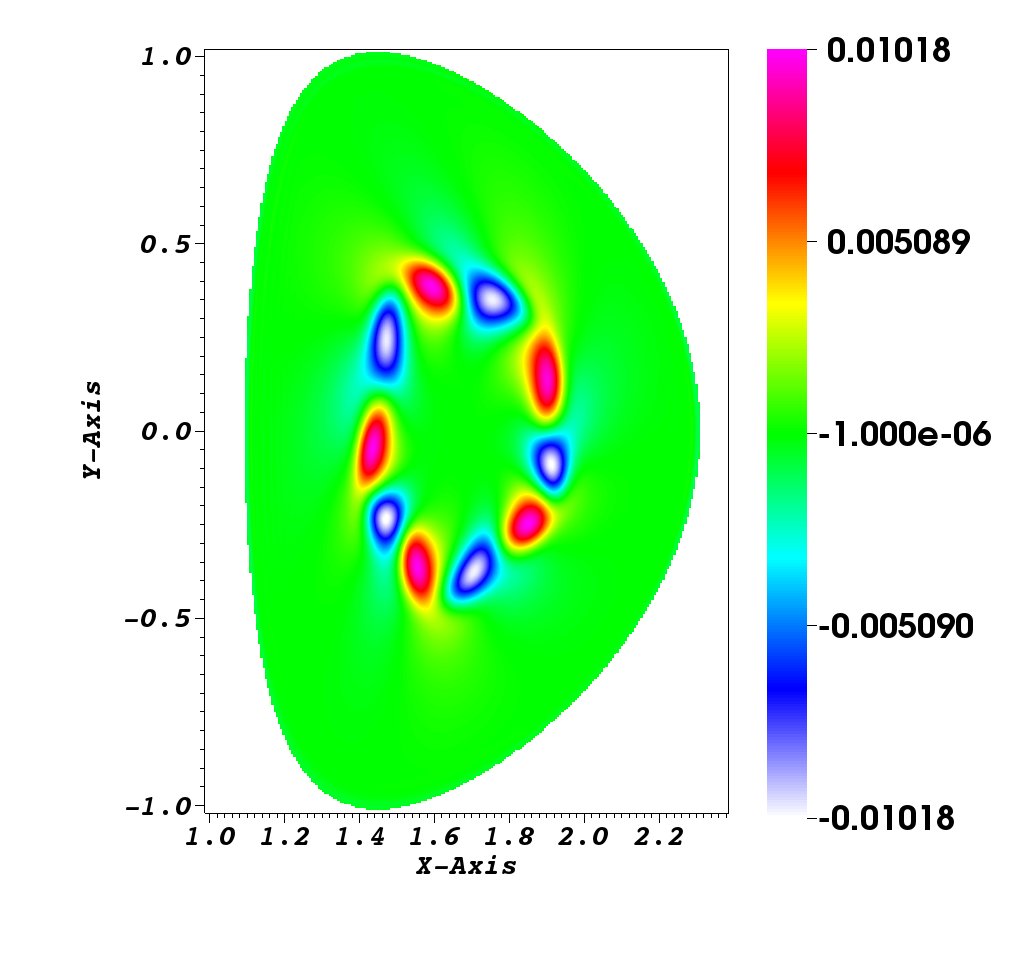} \\
  (a) $t=0$                                      &    (b) $t=100$  \\
  \includegraphics[width=7.5cm]{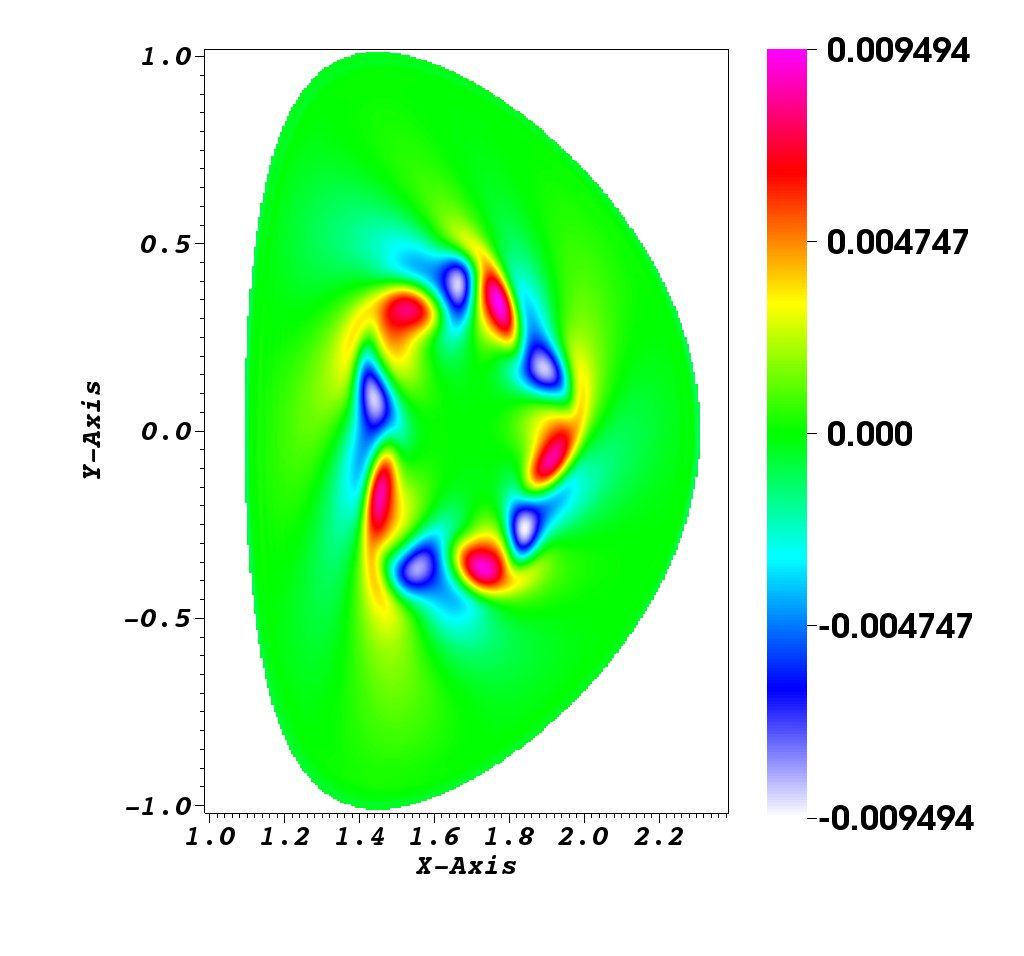} &    \includegraphics[width=7.5cm]{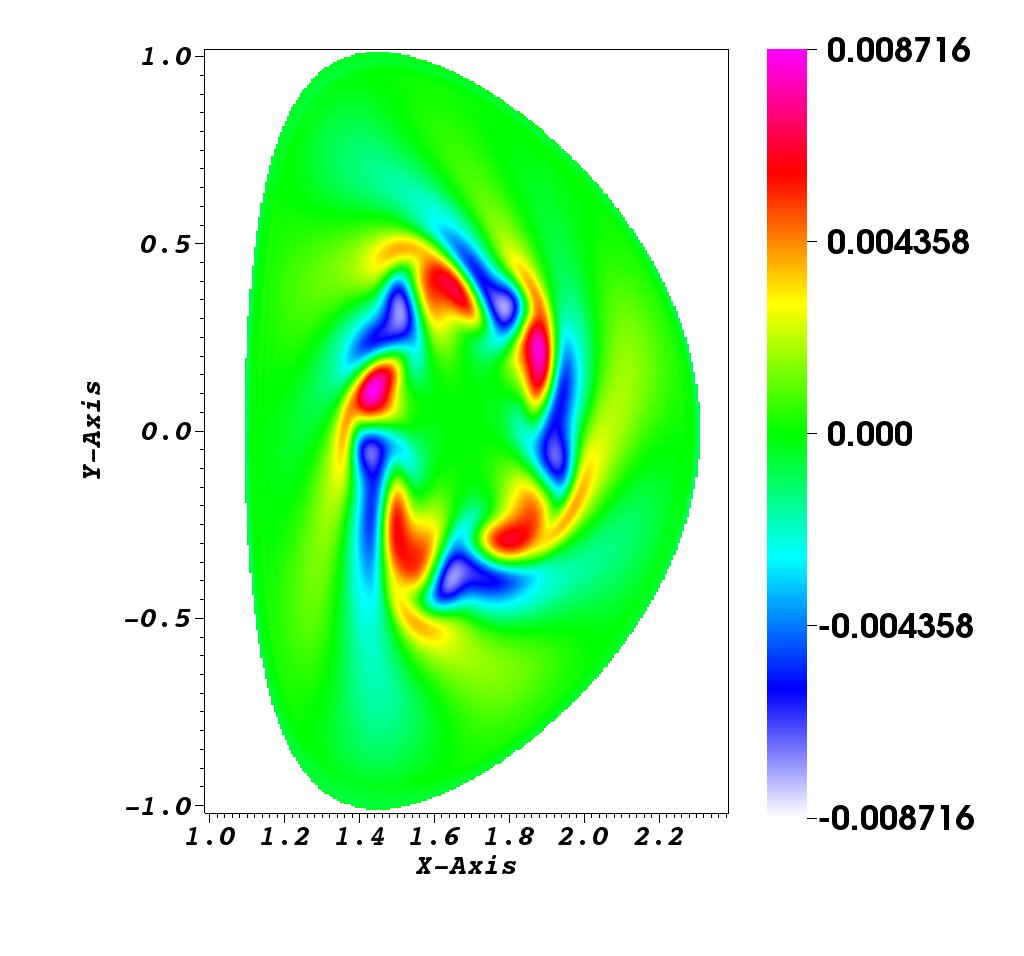} \\
  (c) $t=200$                                   &    (d) $t=300$  
 \end{tabular}
\caption{\label{fig:instability_DK} Instability simulation for guiding-center model in D-shaped domain. The difference between the perturbed density and the steady state density is presented, {\it i.e.} $\delta\rho(t)=\bar\rho(t)-\bar\rho_0$.}
 \end{center}
\end{figure}

Figure~\ref{fig:instability_DK} illustrates the evolution of density governed by the guiding-center model. We present the  difference between the perturbed density and the steady state density, {\it i.e.} $\delta\rho(t)=\bar\rho(t)-\bar\rho_0$. We observe that the difference of density $\delta\rho$ revolves, and  small filaments appear at time $t=200$. Until the time $t=300$, we can  clearly identify the filaments.

\subsection{\bf  Test 3 :  Drift-kinetic model for ion turbulence simulation}
In this subsection, we reproduce the ion turbulence simulation~\cite{bibGV}. This simulation has been realized by different methods but in cylindrical coordinates~\cite{bibGV, bibSamtaney}. We will perform the simulation in Cartesian mesh with the numerical methods presented in section~\ref{sec:schemes}.

The discretization of the Drift-Kinetic model can be developed very similarly as the one for the guiding-center model.
Here, we present some principle discretization steps.

The Vlasov equation of system~\eqref{DK4D} can be split into three equations :
\begin{equation*}
\left\{
\begin{array}{l}
 \displaystyle\frac{\partial {f}}{\partial t}+\mathbf{U}\cdot\nabla_{\mathbf{x}_\bot}{f}=0,\\[3mm]
 \displaystyle\frac{\partial {f}}{\partial t}+  {v}_\|\partial_z f=0,\\[3mm]
 \displaystyle\frac{\partial {f}}{\partial t}+  E_{\|}\partial_{{{v}}_\|}f=0.
 \end{array}
 \right.
\end{equation*}
Thus when applying the Semi-Lagrangian method, we use the Strang splitting method~\cite{bibStrang} for time discretization, and a second order predictor-corrector method for searching the characteristic foot. The Semi-Lagrangian method is suitable for linear phase simulation, because it doesn't have CFL constraint and it is very accurate. However, it doesn't preserve well the conservation properties of  physical models~\cite{YF15} in nonlinear phase. Therefore, we should use the conservative finite difference method presented in section~\ref{sec:HFD}, where the 4th order Runge-Kutta method is used for time discretization.  The criterion to pass from the Semi-Lagrangian to the finite difference methods is as follows
\begin{equation}
 \left|\int_{\Omega}\int_{\mathbb{R}}\left[f(t_n)-f(t_{n-1})\right]d v d\xx \right|\,>\, h^3,
\label{eq:criterion}
\end{equation}
where $h$ is the smallest space step. 

The  quasi-neutrality equation of the system~\eqref{DK4D} is a three-dimensional elliptic problem.
Obviously, a direct resolution of this equation is very costly.
However, we notice that the diffusion term $\frac{\rho _0(\mathbf{x}_\bot)}{B}$ only depends  on $\mathbf{x}_\bot$,
and  the second term of  the quasi-neutrality equation is vanishing by taking average  in $z$-direction.
Thus  averaging the the quasi-neutrality equation  in $z$-direction, 
we get a 2D average equation
\begin{equation}
 -\nabla_\bot\cdot\left(\frac{\rho _0(\mathbf{x}_\bot)}{B}\nabla_\bot\bar{\phi}\right)=\bar{\rho}-\rho_0 \quad\text{ in }D.
 \label{eq:average_admin}
\end{equation}
Taking difference between the quasi-neutrality equation and the average equation, it yields a fluctuation equation :
\begin{equation}
 -\nabla_\bot\cdot\left(\frac{\rho _0(\mathbf{x}_\bot)}{B}\nabla_\bot{\phi'}\right)+\frac{\rho_0(\mathbf{x}_\bot)}{T_e(\mathbf{x}_\bot)}\phi'=\rho'-\bar{\rho} \quad\text{ in } \Omega=D\times[0, L_z],
  \label{eq:fluctuation_admin}
\end{equation}
with  $\phi'=\phi-\bar{\phi}$.
Note that the terms $\frac{\rho _0(\mathbf{x}_\bot)}{B}$ and $\frac{\rho_0(\mathbf{x}_\bot)}{T_e(\mathbf{x}_\bot)}$ are independent on  $z$. 
Thus the fluctuation equation~\eqref{eq:fluctuation_admin} can be solve slice by slice in $z$-direction.
Moreover, the Dirichlet boundary conditions can also be decomposed for the average equation~\eqref{eq:average_admin} and the fluctuation equation~\eqref{eq:fluctuation_admin} as follows
\begin{equation*}
\left\{
\begin{array}{ll}
 \bar{\phi}=0,&\,\forall\mathbf{x}_\bot\in\partial D,\\[3mm]
 \phi'=0,&\,\forall\mathbf{x}_\bot\in\partial D\times[0,L_z].
 \end{array}
 \right.
\end{equation*}
Therefore, the resolution for Poisson equation in Section~\ref{sec:poisson} can be applied for~\eqref{eq:average_admin} and~\eqref{eq:fluctuation_admin}.

Now we perform the ion turbulence simulation. The plasma is initialized by exciting a single ion temperature gradient (ITG) model $(m,n)$ (where $m$ is a poloidal mode and $n$ is a toroidal mode).
The distribution function is thus considered at the initial time as the sum of an equilibrium and a perturbed part: $f=f_{\text{eq}}+\delta f$.
The equilibrium part $f_{\text{eq}}$ is chosen as a local Maxwellian
\begin{equation*}
 f_{\text{eq}}(r,v_\|)=\frac{n_0(r)}{(2\pi T_i(r))^{1/2}}\exp\left(-\frac{v^2_\|}{2T_i(r)}\right),
\end{equation*}
while the perturbation $\delta f$ is determined as
\begin{equation*}
 \delta f=f_{\text{eq}}\, \varepsilon\exp\left( - \frac{(r - r_p)^2}{\delta r}  \right) \, \cos \left( \frac{2\pi n}{L} z + m \theta \right),
\end{equation*}
where the profiles $n_0(r)$, $T_i(r)$ and $T_e(r)$  satisfy 
\begin{equation*}
  \frac{\partial_r P(r)}{P(r)}\, = \, -\kappa_P \cosh^{-2} \left( \frac{r - r_p}{\delta r_P} \right),\, \text{for }P = n_0,\, T_i \text{ and } T_e,
\end{equation*}
together with the normalization
\begin{equation*}
  \int_{r_{\min}}^{r_{\max}}  n_0 (r)\, dr = r_{\max} - r_{\min}, \, T_i(r_p) = T_e(r_p) = 1.
\end{equation*}
This gives the formulas
\begin{equation*}
  P(r) = C_P \exp \left( -\kappa_P\, \delta r_P \tanh\left( \frac{r - r_p}{\delta r_P} \right)  \right),
\end{equation*}
where $C_{T_i} = C_{T_e} = 1$ and $C_{n_0} = \frac{r_{\max} - r_{\min}}{\int_{r_{\min}}^{r_{\max}} \exp\left( -\kappa_{n_0} \delta r_{n_0} \tanh\left( \frac{r - r_p}{\delta r_{n_0}} \right) \right) dr }$.

In this simulation, we choose the following parameters
\begin{equation*}
  \begin{array}{c}
  r_{\min} = 0,\, r_{\max} = 14.5, \, \kappa_{n_0} = 0.055, \kappa_{T_i} = \kappa_{T_e} = 0.27586,\\[3mm]
  \delta r_{T_i} = \delta r_{T_e} = \frac{\delta r_{n_0}}{2} = 1.45,\, \varepsilon = 10^{-6},\\[3mm]
  n = 1,\, m=5, \, L=1506.759067, \, v_{\max} = 8, \, r_p = \frac{r_{\max} + r_{\min}}{2},\, \delta r = \frac{4\delta r_{n_0}}{\delta r_{T_i}}.
\end{array}
\end{equation*}

Let us first compare the different discretization methods. The Semi-Lagrangian methods with cubic Hermite reconstruction and the HWENO reconstruction are used to solve the 4D Drift-Kinetic model, with small time step such that the CFL number is small than 1. Then we compare the numerical results with the one obtained by the mixed Semi-Lagrangian/finite difference method, where large time step (CFL$>1$) is used for the Semi-Lagrangian method in linear phase and small time step  (CFL$<1$) is used in nonlinear phase. We emphasize that the Semi-Lagrangian method switches to the finite difference method automatically by the criterion~\eqref{eq:criterion}.

\begin{figure}
 \begin{center}
  \begin{tabular}{cc}
   \includegraphics[width=7cm]{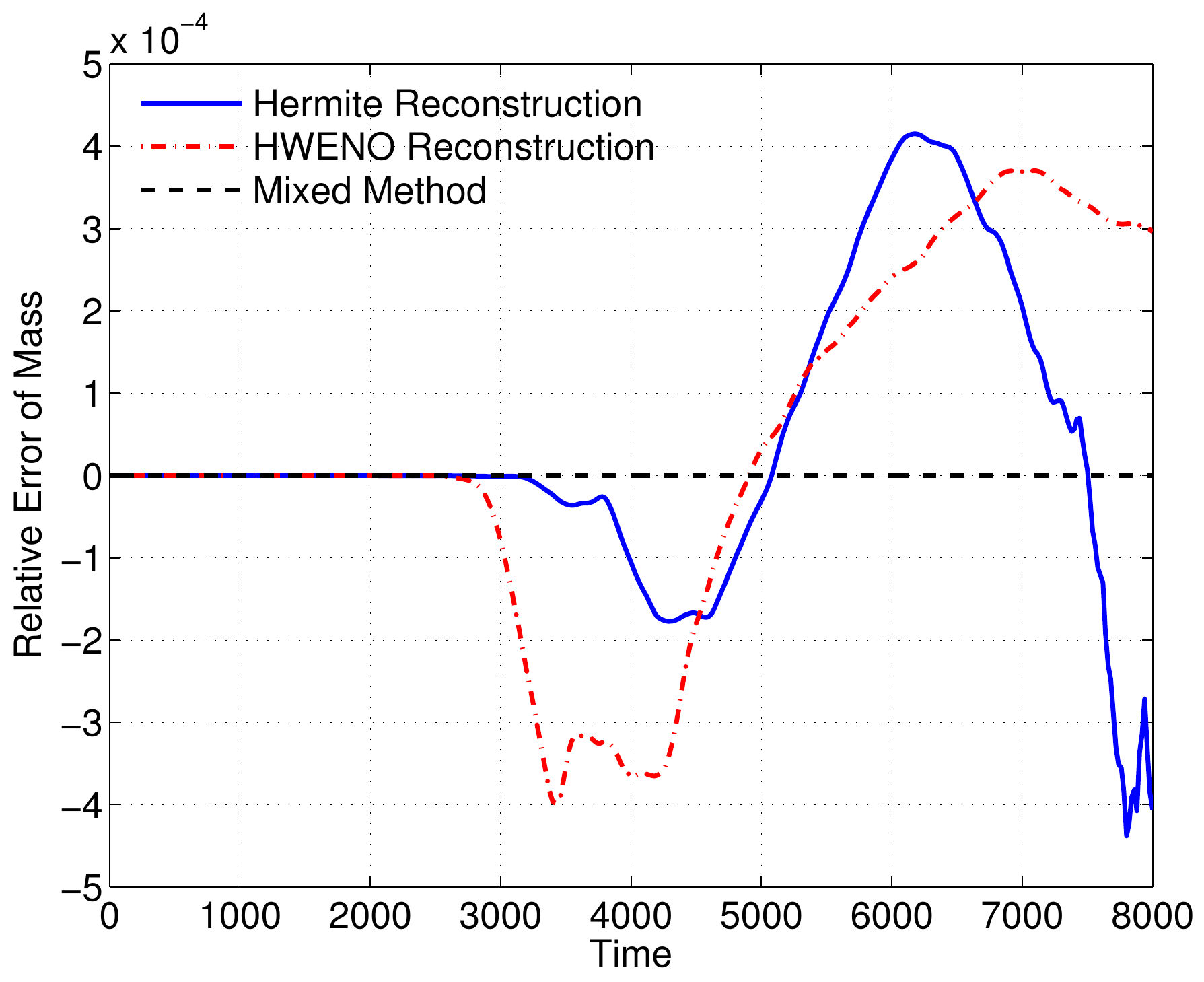}     &    \includegraphics[width=7cm]{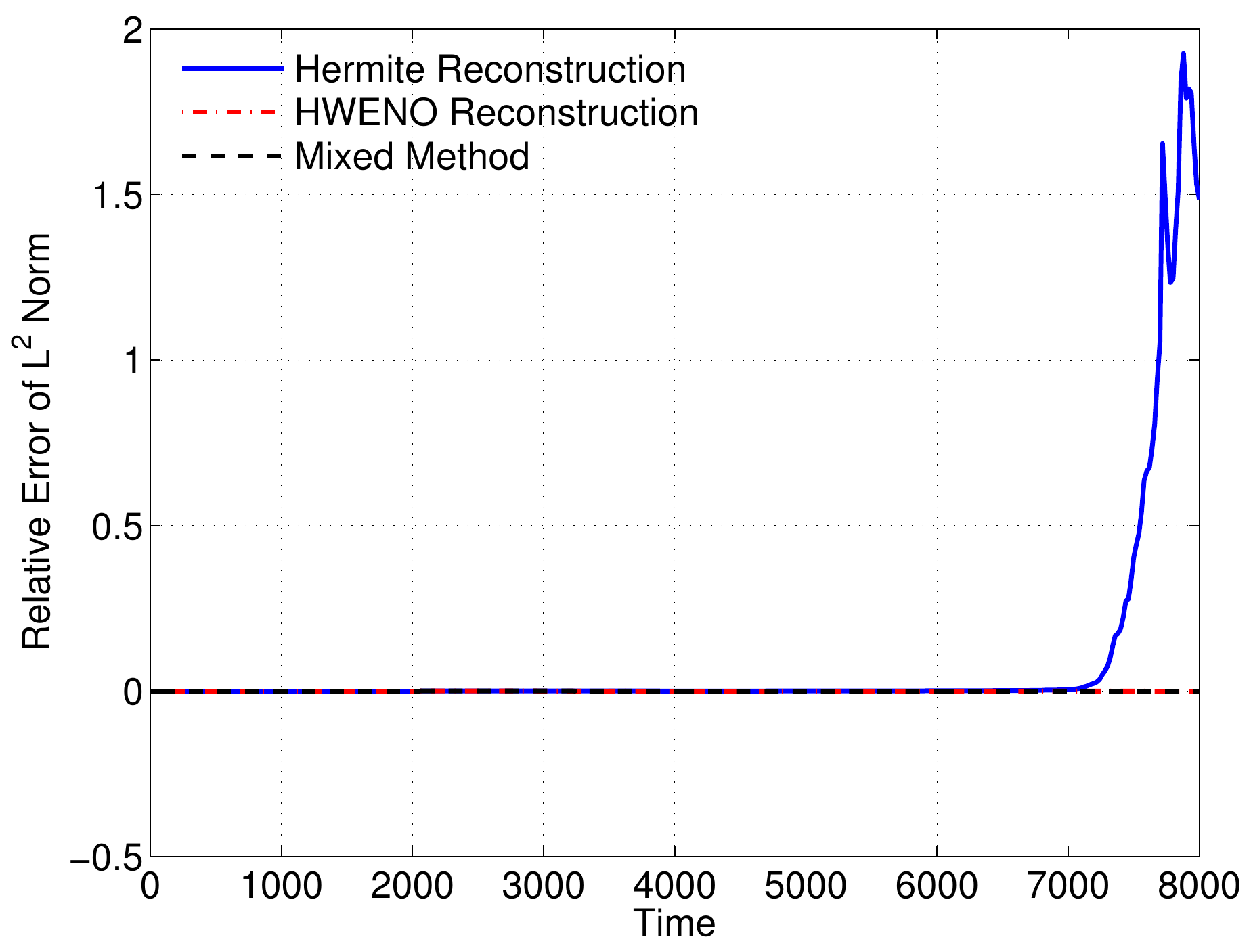} \\
    (a)                                                &    (b)  \\
    \includegraphics[width=7cm]{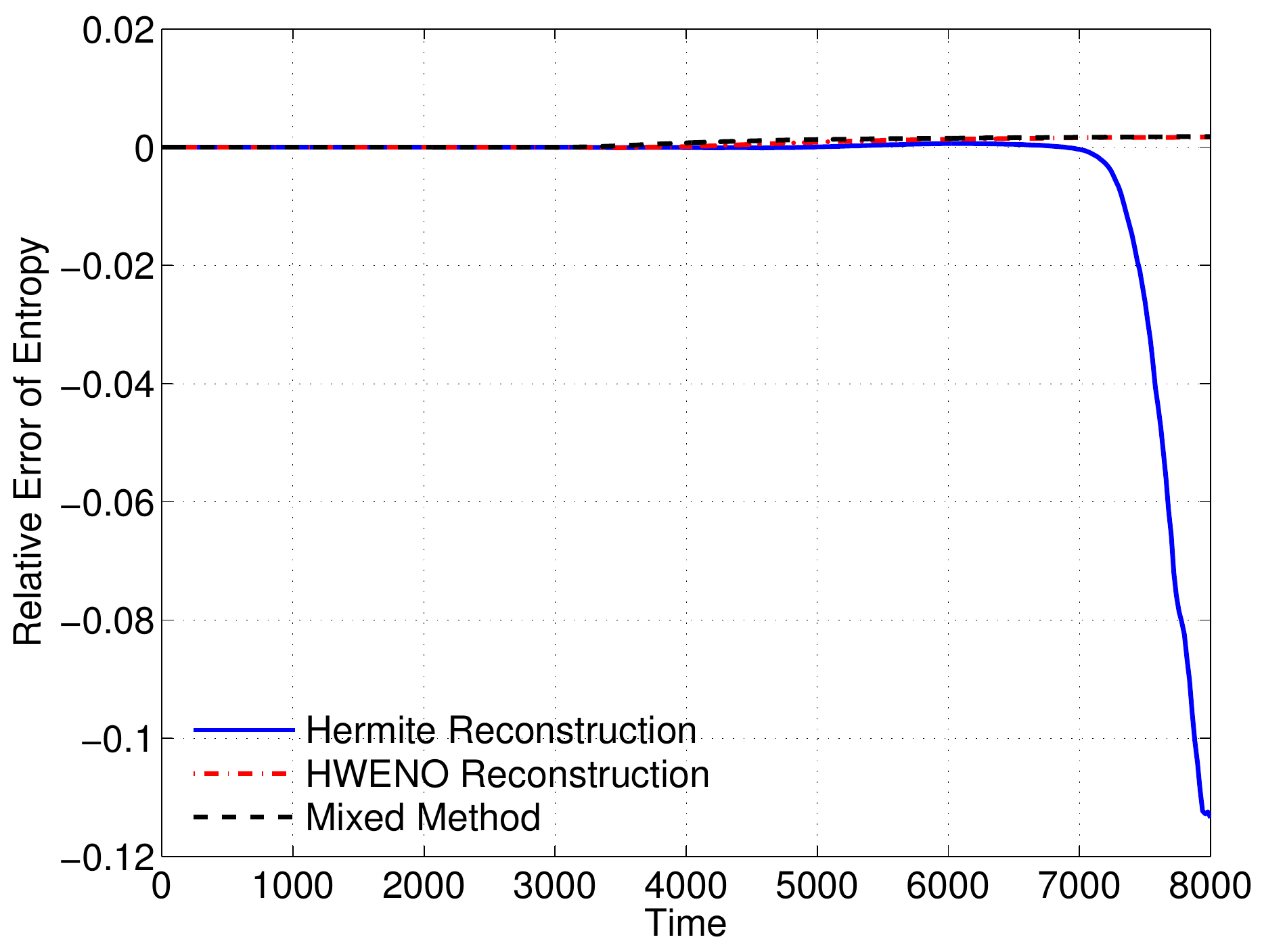} &   \includegraphics[width=7cm]{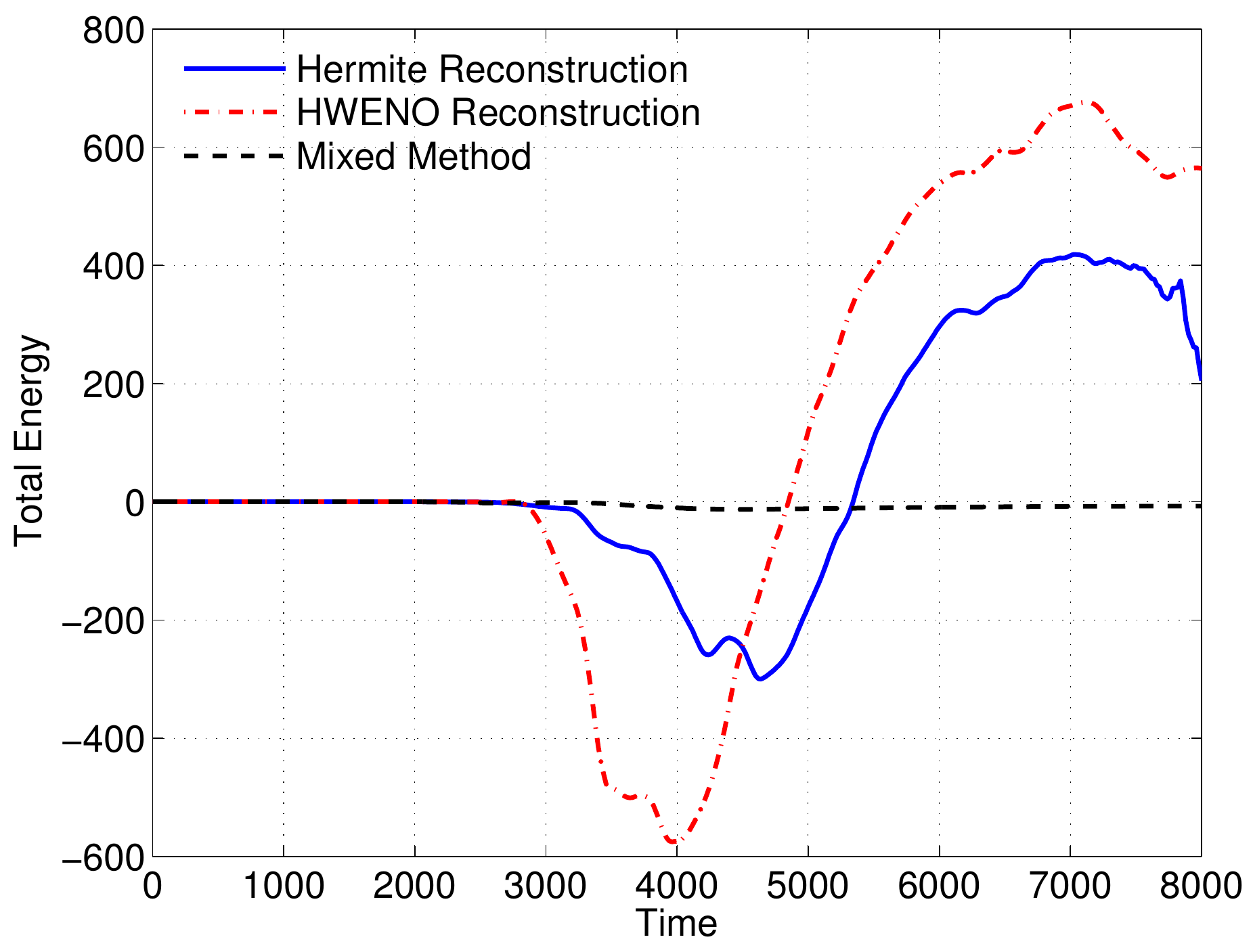}  \\
    (c)                                                &    (d)
  \end{tabular}
   \caption{\label{fig:DK_compare}Comparison of different reconstruction methods for the conservation laws for the Drift-Kinetic model. The mesh size is $n_x\times n_y\times n_z\times n_v = 128\times 128\times 32 \times 65$.}
 \end{center}
\end{figure}

In Figure~\ref{fig:DK_compare}, we summarize relative errors of the conservation laws for the Drift-Kinetic model for the different methods. We notice that these three methods have almost the same results when $t<3000$, while these results differ significantly when $t\geq 3000$. We thus denote the linear phase for  $t<3000$ and the nonlinear phase for  $t\geq 3000$.

 The Semi-Lagrangian methods can not conserve well the mass in the nonlinear phase, while the finite difference method conserve exactly the mass (see Figure~\ref{fig:DK_compare}(a)). Then from Figures~\ref{fig:DK_compare}(b),~\ref{fig:DK_compare}(c), we observe that the Semi-Lagrangian method with Hermite reconstruction loses completely the conservation properties for $L^2$ norm and entropy for long time simulation, since it involves too much spurious oscillation. At contrast, the Semi-Lagrangian method and the finite difference method with HWENO reconstruction work much better. Finally, we see  the  mixed Semi-Lagrangian/finite difference method has better energy conservation property than the Semi-Lagrangian methods in the nonlinear phase (see Figure~\ref{fig:DK_compare}(d)). Therefore,  the  mixed Semi-Lagrangian/finite difference method is better than the  Semi-Lagrangian methods for long time ion turbulence simulation. Moreover, since the  Semi-Lagrangian method is used in linear phase, thus our mixed method is more efficient than the pure finite difference method.

\begin{figure}
 \begin{center}
  \begin{tabular}{cc}
   \includegraphics[width=7cm]{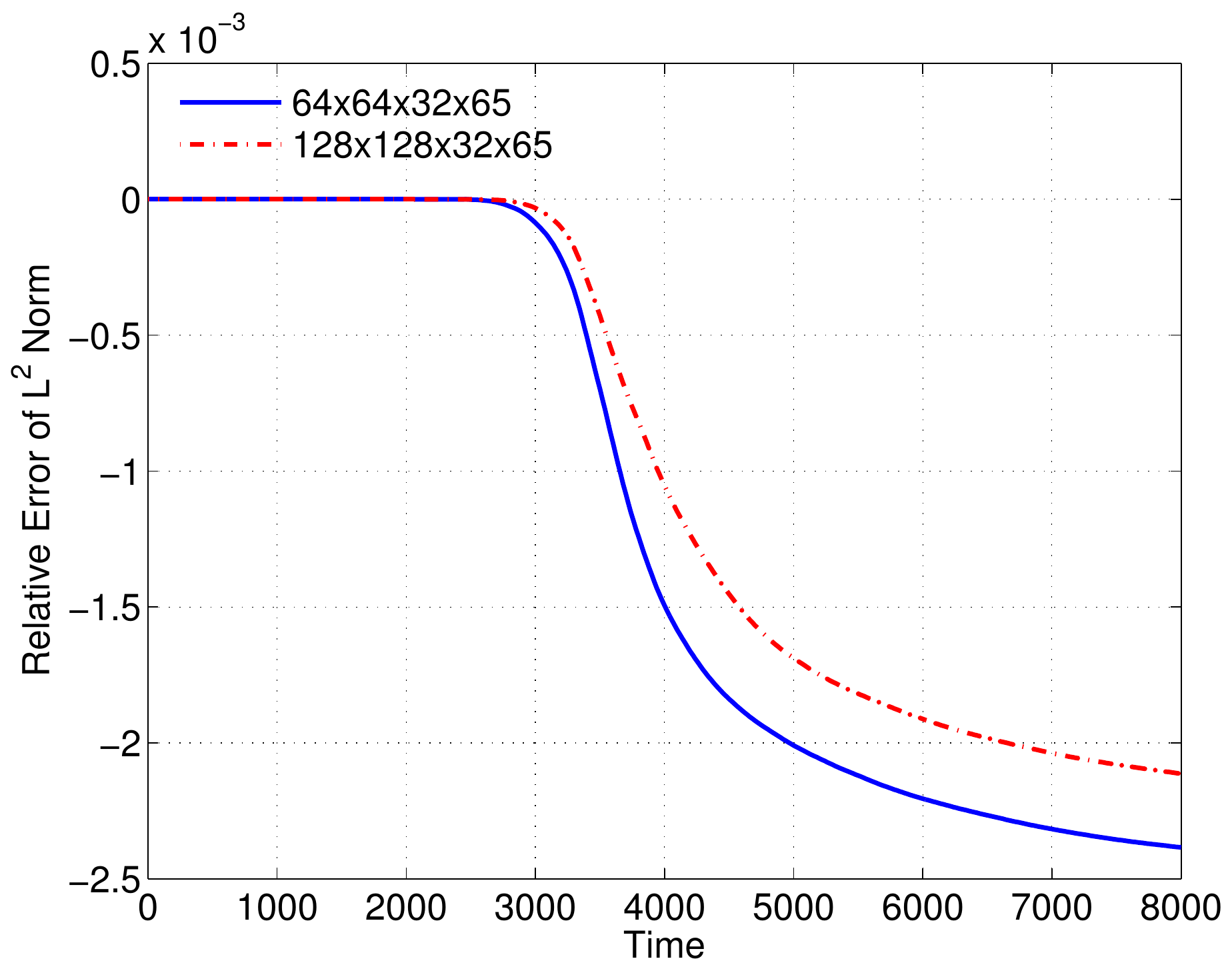} &    \includegraphics[width=7cm]{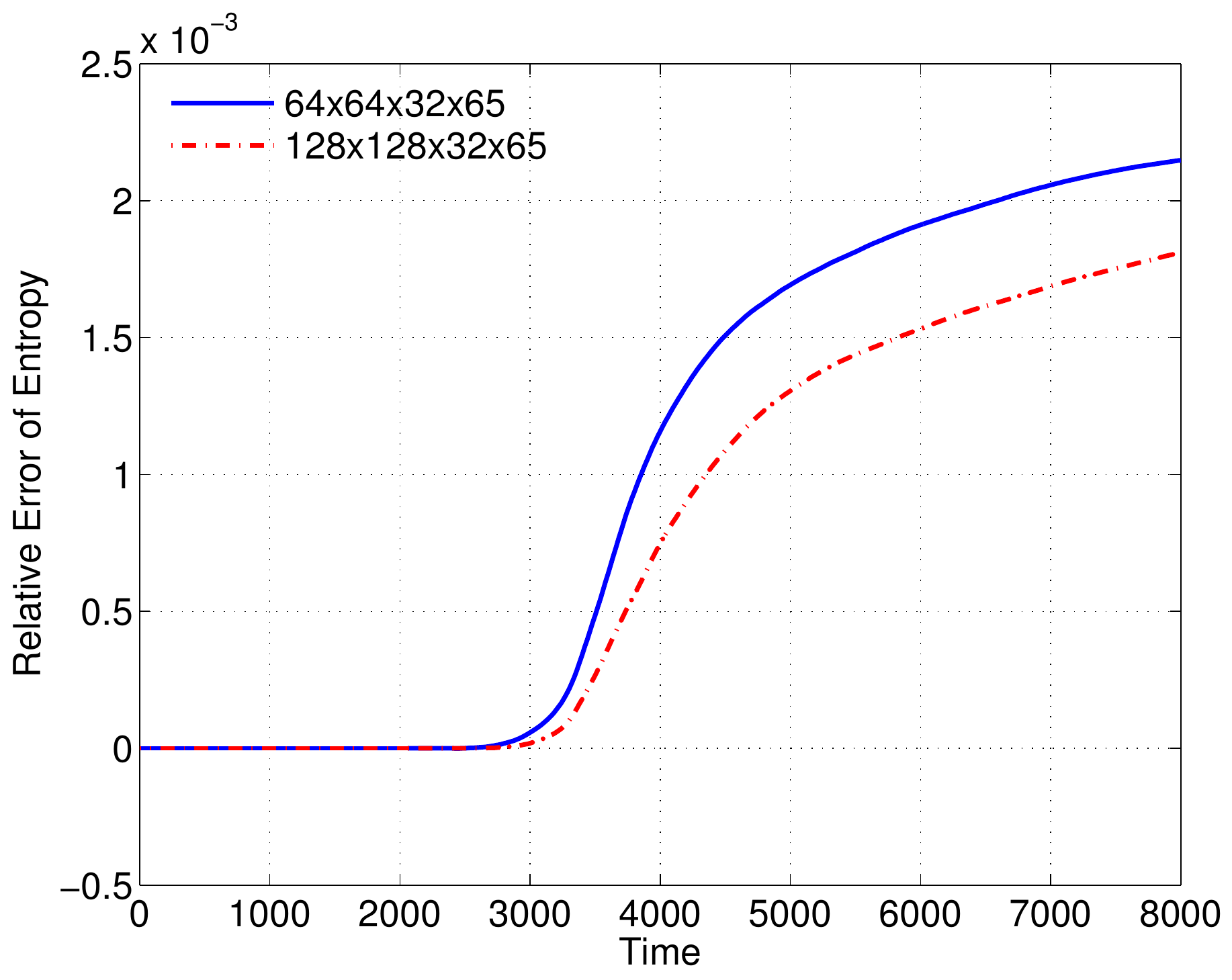} \\
    (a)                                          &    (b)  \\
   \includegraphics[width=7cm]{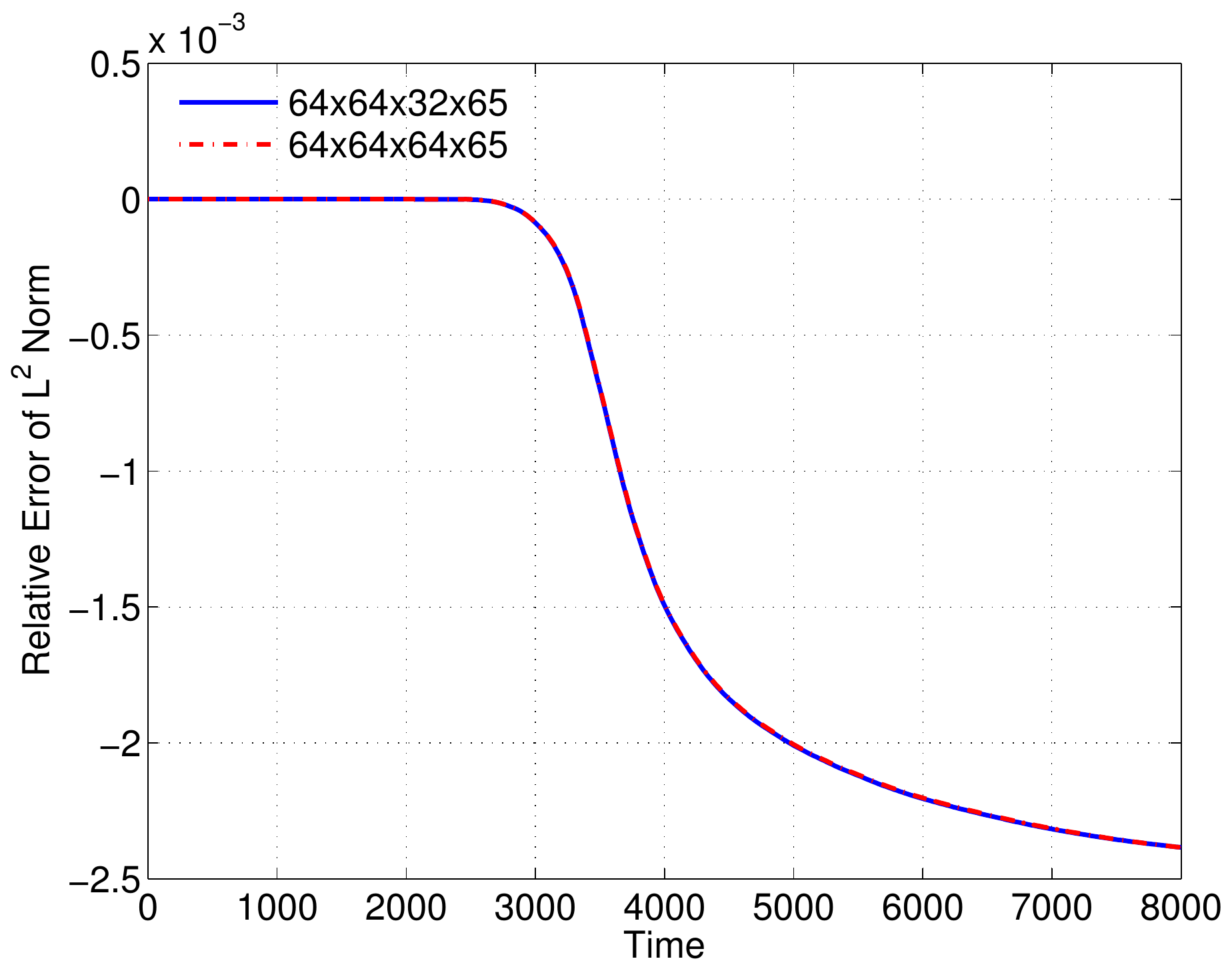} &    \includegraphics[width=7cm]{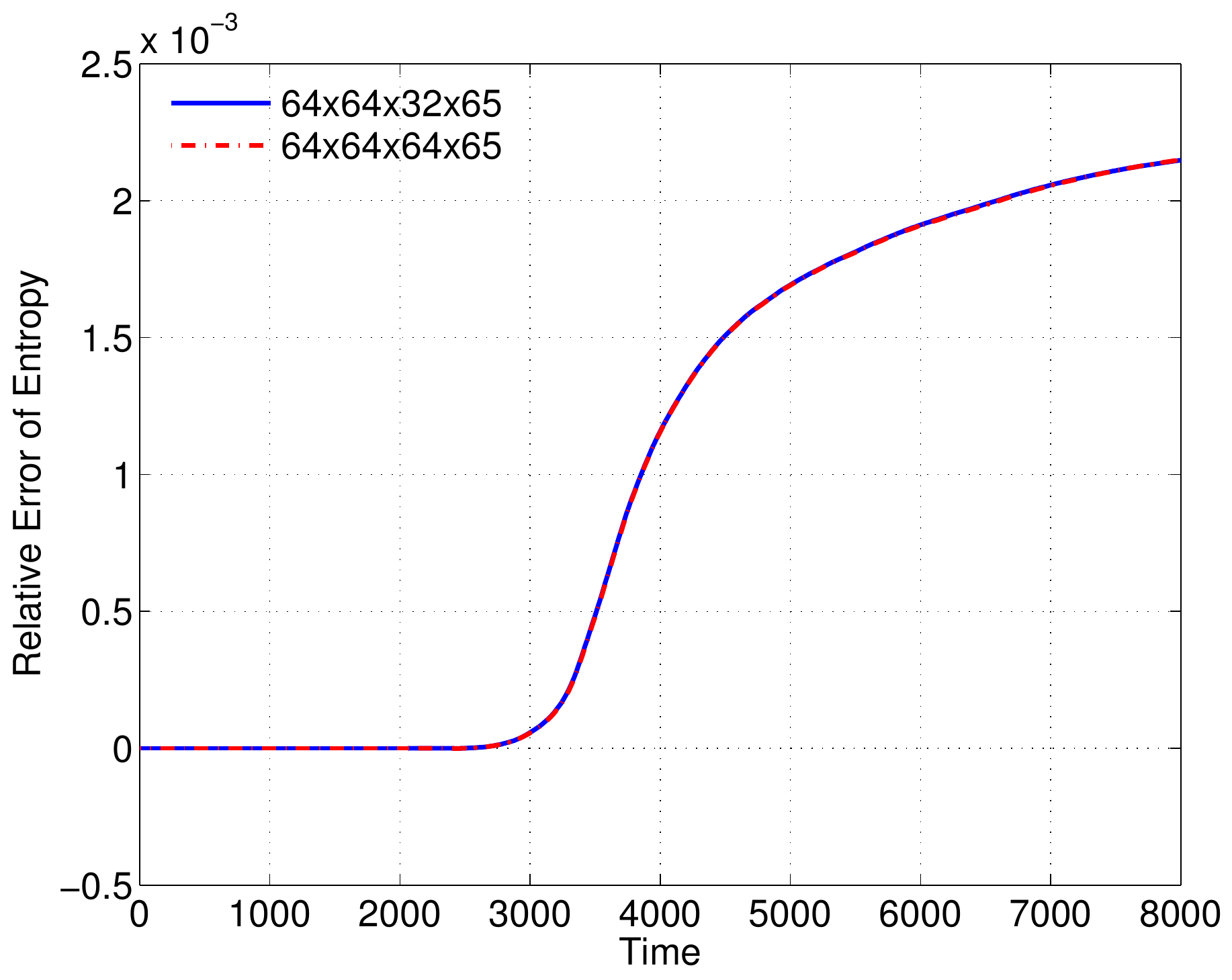} \\
    (c)                                          &    (d)  \\
   \includegraphics[width=7cm]{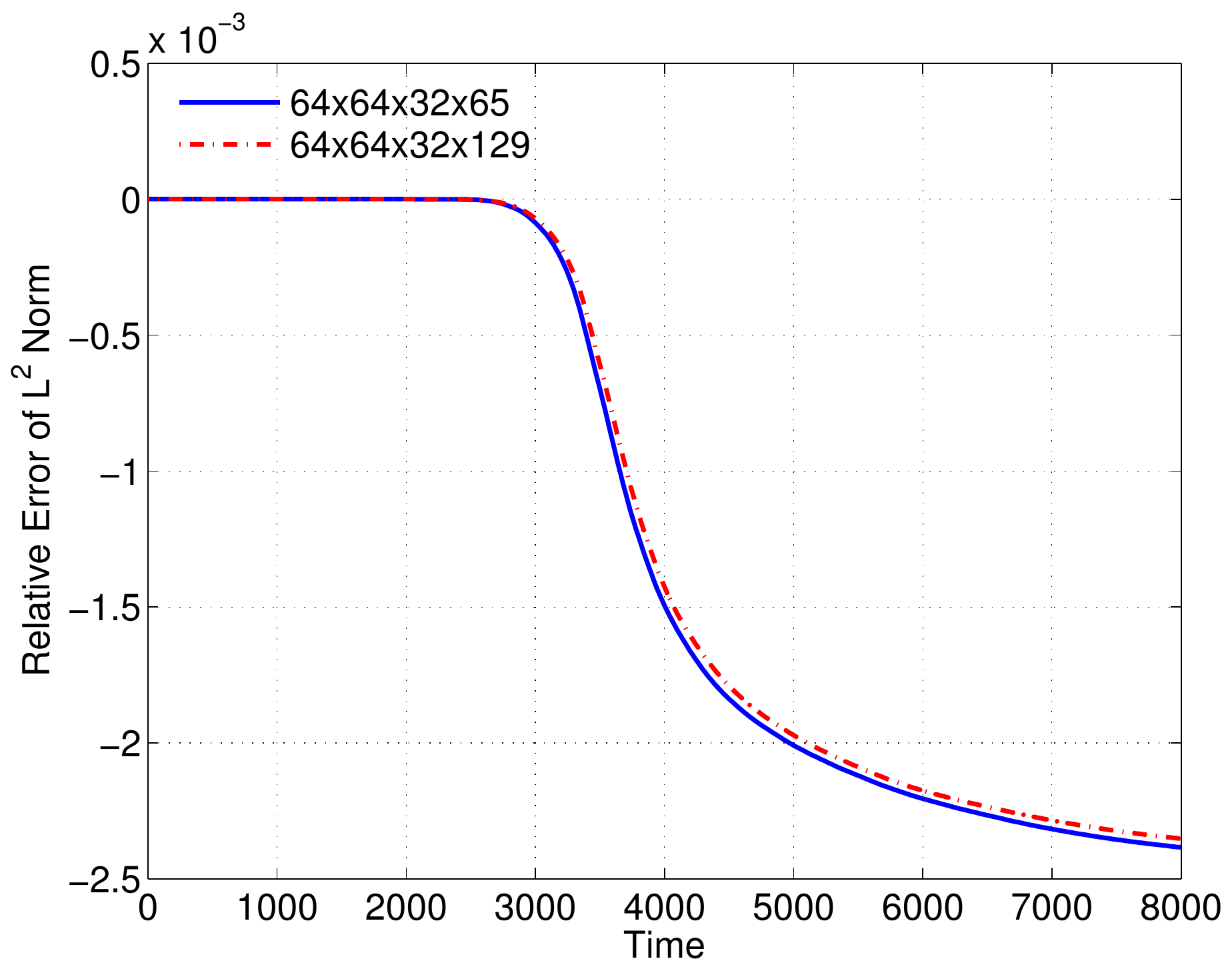} &    \includegraphics[width=7cm]{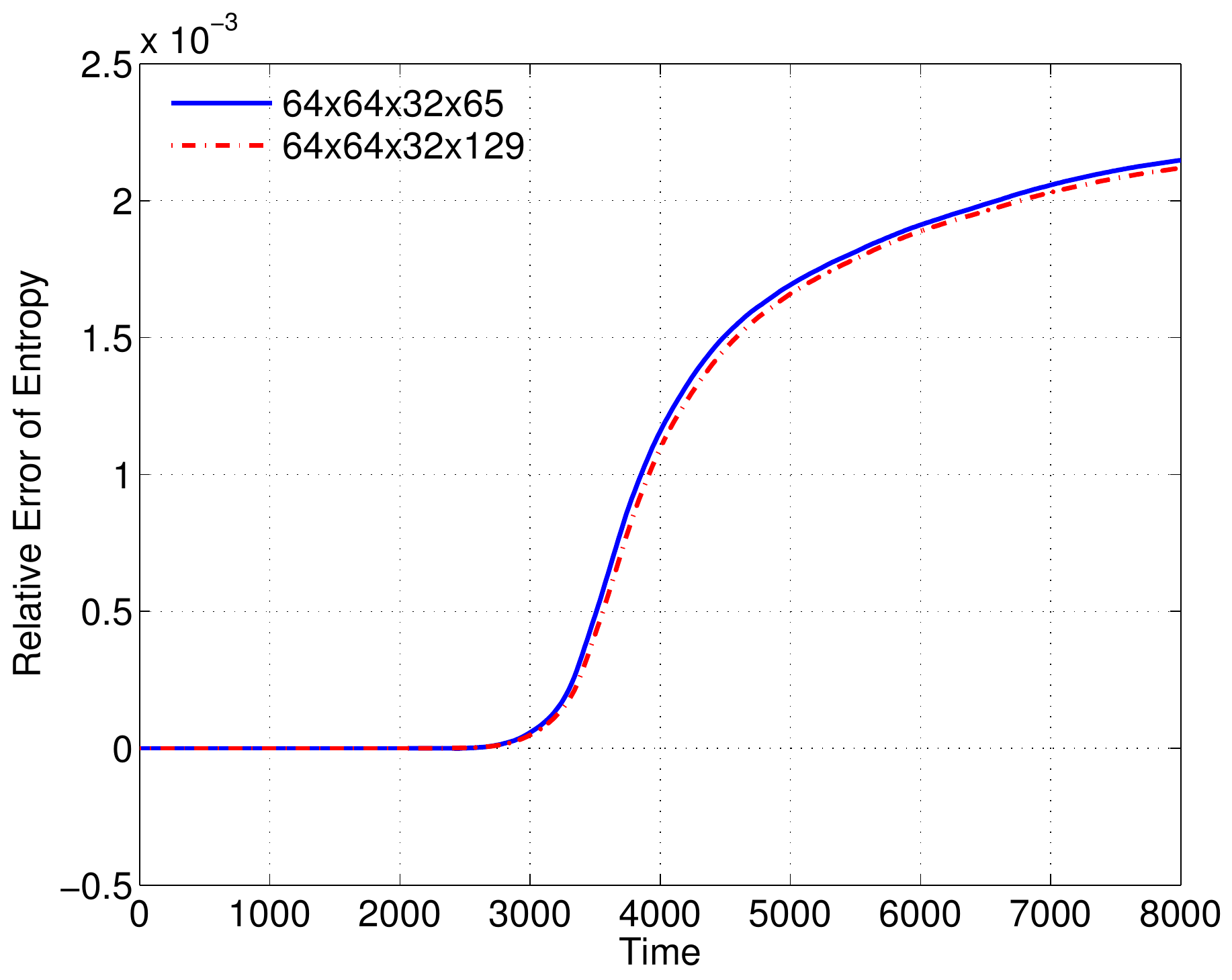} \\
    (e)                                          &    (f)  
  \end{tabular}
   \caption{\label{fig:DK_refine}Refinement of mesh size in different directions for  $L^2$ norm and entropy conservations of the Drift-Kinetic model. Mixed Semi-Lagrangian/finite difference method is used.}
 \end{center}
\end{figure}

\begin{figure}
 \begin{center}
  \begin{tabular}{cc}
    \includegraphics[width=7cm]{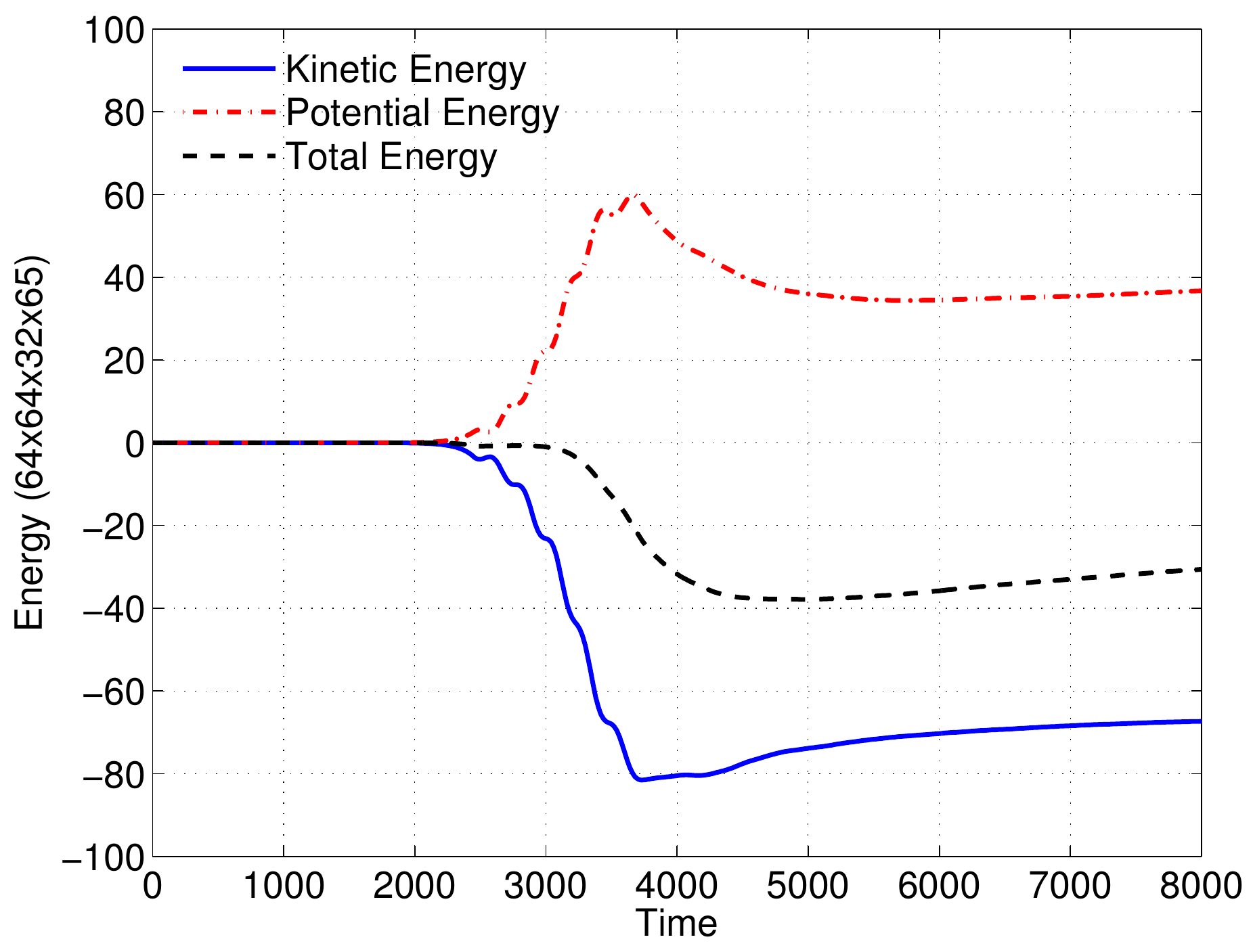} &    \includegraphics[width=7cm]{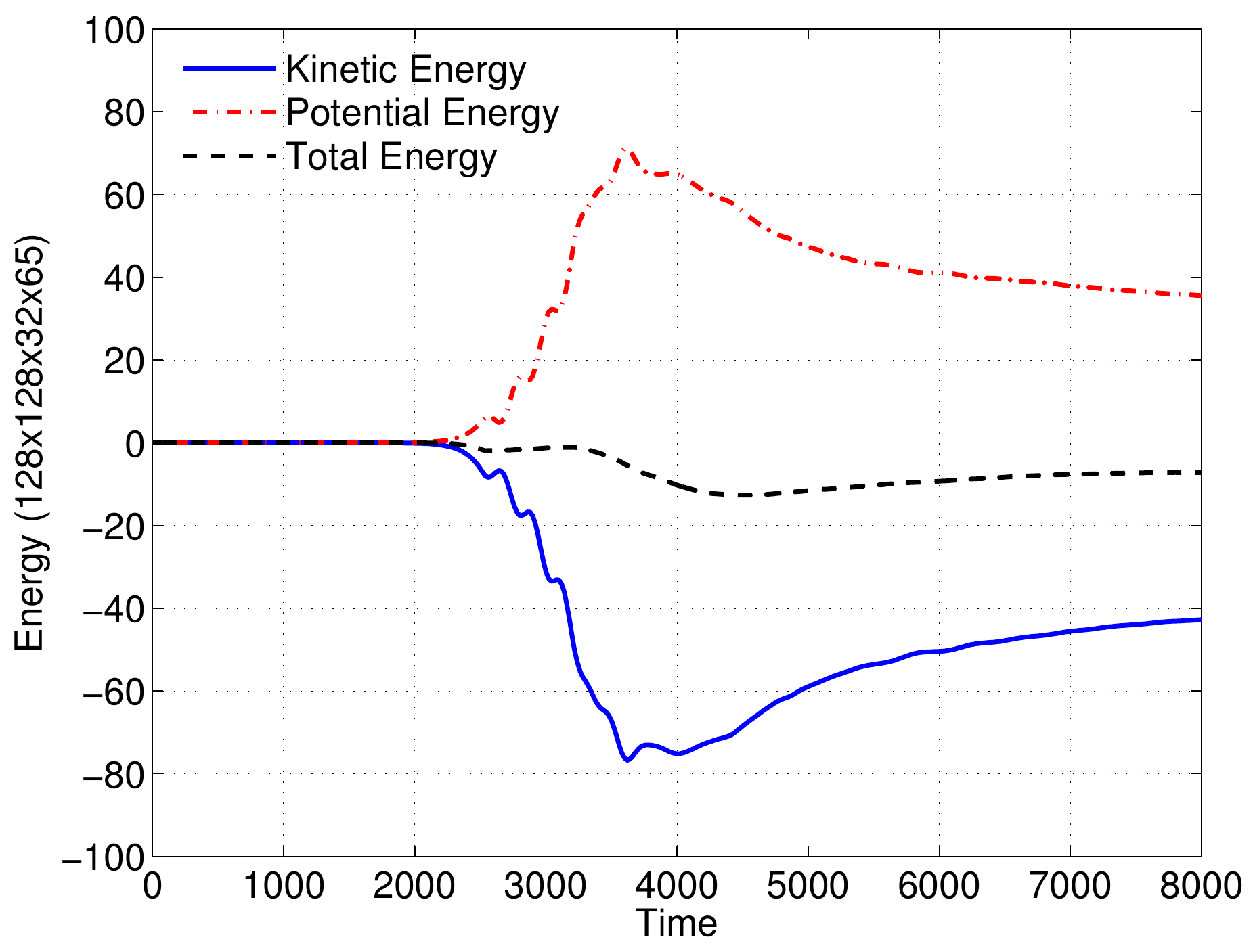} \\
    (a)                            &    (b)  \\
    \includegraphics[width=7cm]{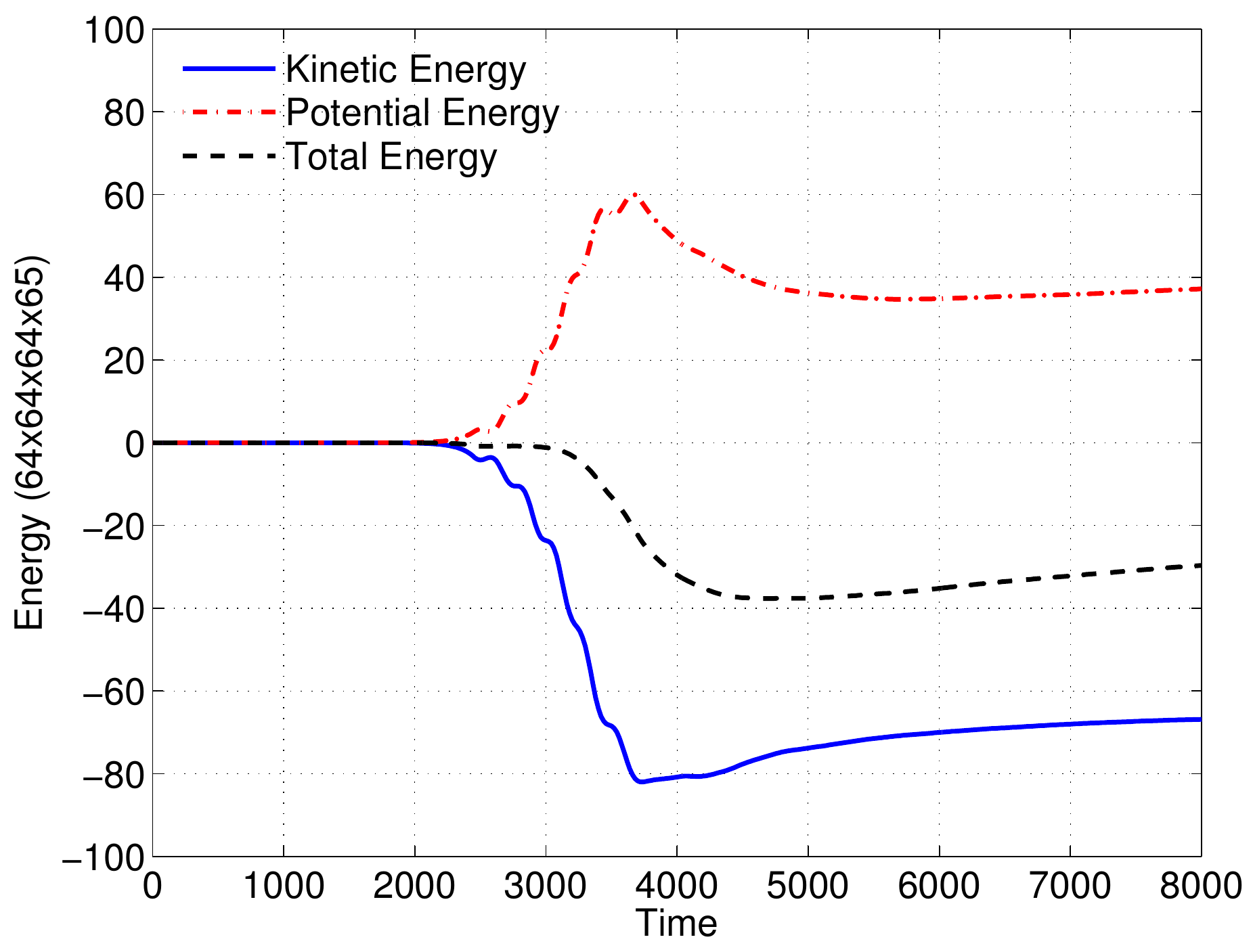} &    \includegraphics[width=7cm]{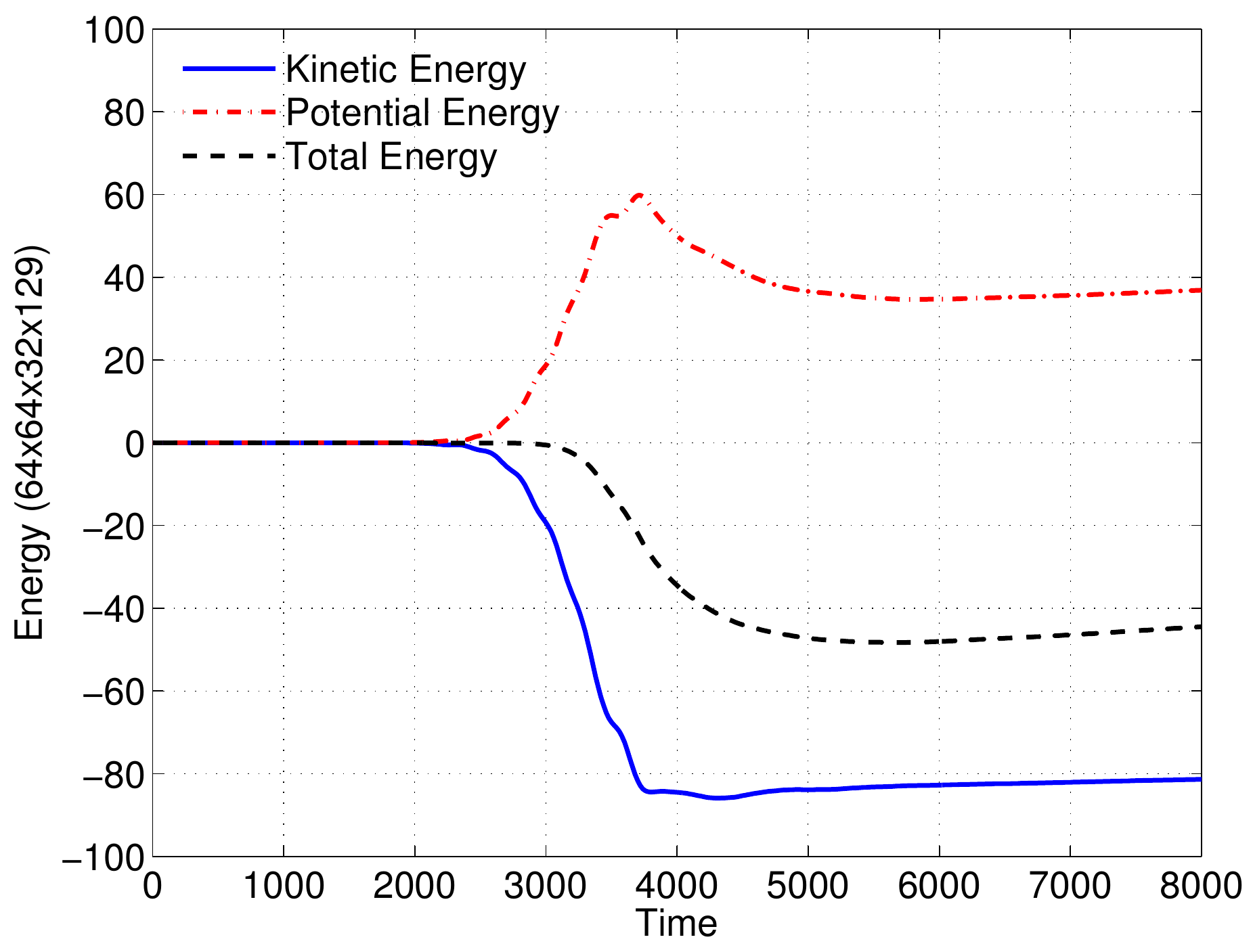} \\
    (c)                             &    (d)
  \end{tabular}
   \caption{\label{fig:DK_refine_energy}Refinement of mesh size in different directions for  the energy conservations of the Drift-Kinetic model. Mixed Semi-Lagrangian/finite difference method is used.}
 \end{center}
\end{figure}

We next investigate the $L^2$ norm and entropy convergence properties of the mixed Semi-Lagrangian/finite difference method by refining mesh size in different directions. The results calculated with mesh size  $n_x\times n_y\times n_z\times n_v = 64\times 64\times 32 \times 65$ is used as a reference solution.  In Figure~\ref{fig:DK_refine}, we observe that a significant improvement is obtained by refining in $x,\,y$ directions. The refinement in $z$ direction doesn't improve the results.  The ones obtained by refining in $v$ direction is slightly better than the reference solution. Figure~\ref{fig:DK_refine_energy} presents the refinement results for energy conservations. We see again the refinement  in $x,\,y$ directions improves best the energy conservations.

\begin{figure}
 \begin{center}
   \includegraphics[width=7cm]{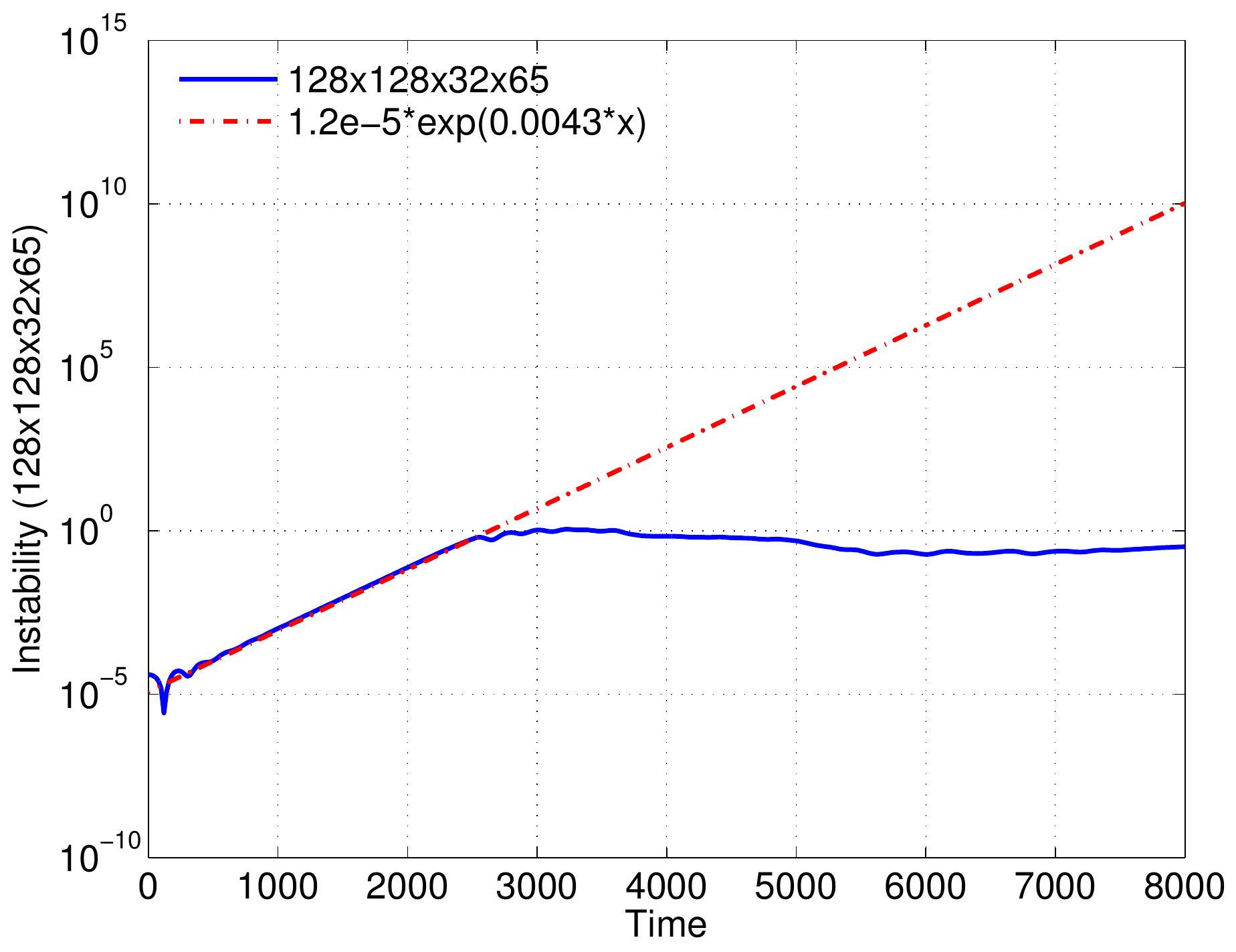} 
   \caption{\label{fig:DK_instability}Growth rate of instability. Mixed Semi-Lagrangian/finite difference method is used.}
 \end{center}
\end{figure}

At last, we present the evolution of the distribution function during ion turbulence simulation. We first notice in Figure~\ref{fig:DK_instability} that the instability develops exponentially in the linear phase, where the growth of instability is measured by the quantity
\begin{equation*}
\sqrt{\int\phi(t,x,y,z)^2dx\,dy\,dz},
\end{equation*}
where $(x,y)\in\{(x,y)\,:\,x^2+y^2 = r_p^2\}$. Then the instability reaches a saturation point, which corresponds to the starting point of the nonlinear phase. 

\begin{figure}
\begin{center}
 \begin{tabular}{cc}
  \includegraphics[width=6.5cm]{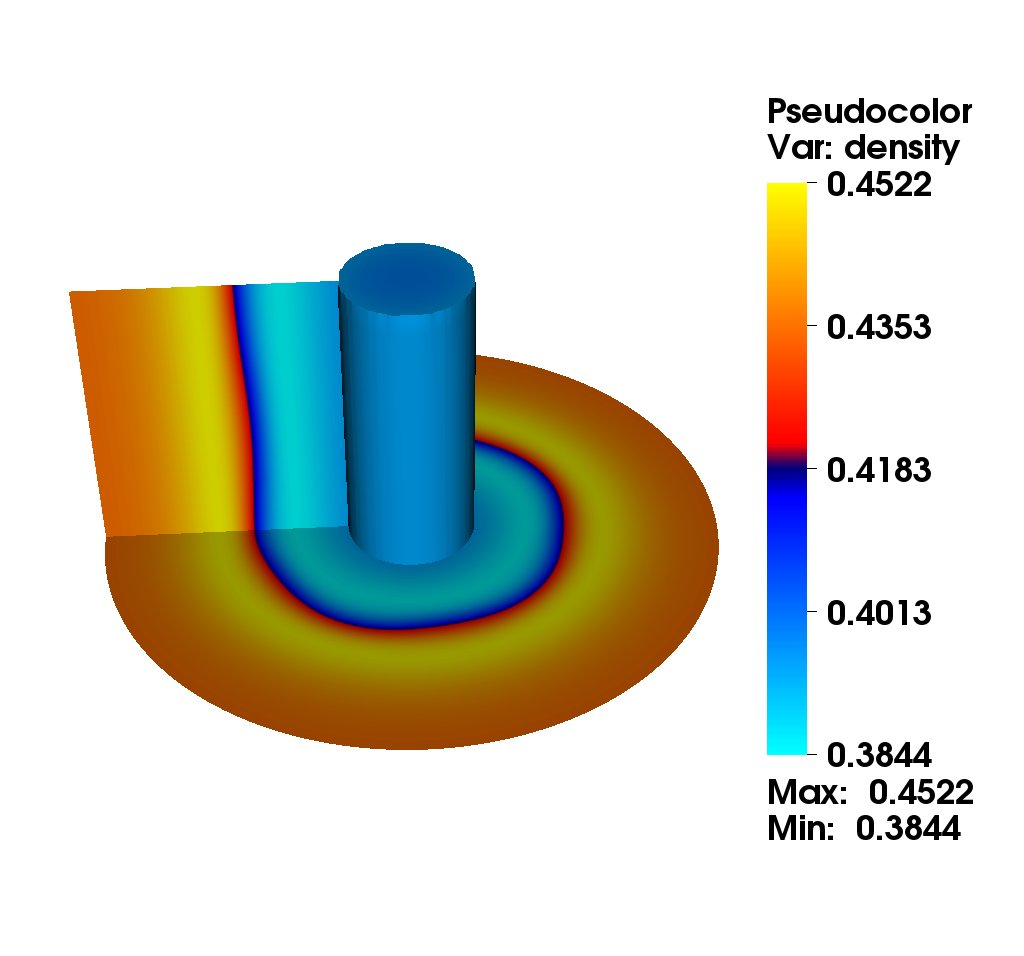} &    \includegraphics[width=6.5cm]{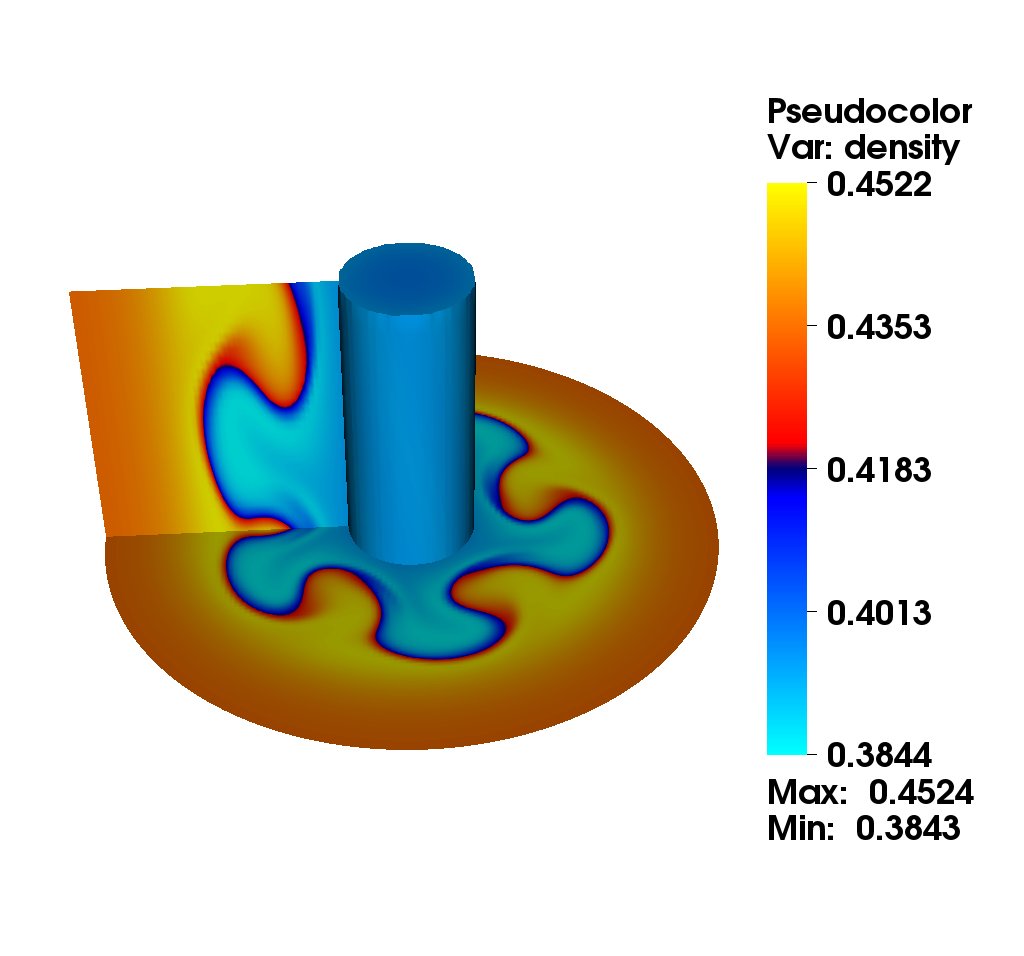} \\
  (a) $t=2000$                                      &    (b) $t=3000$  \\
  \includegraphics[width=6.5cm]{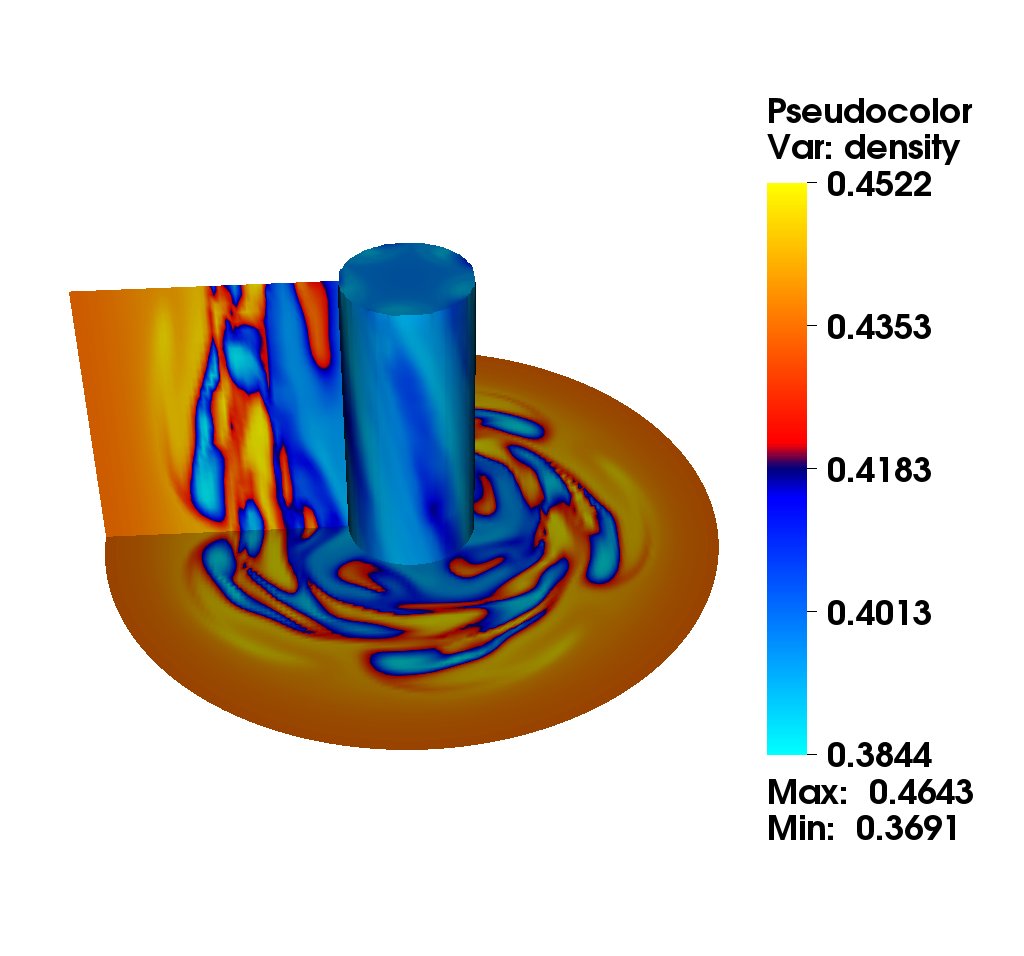} &    \includegraphics[width=6.5cm]{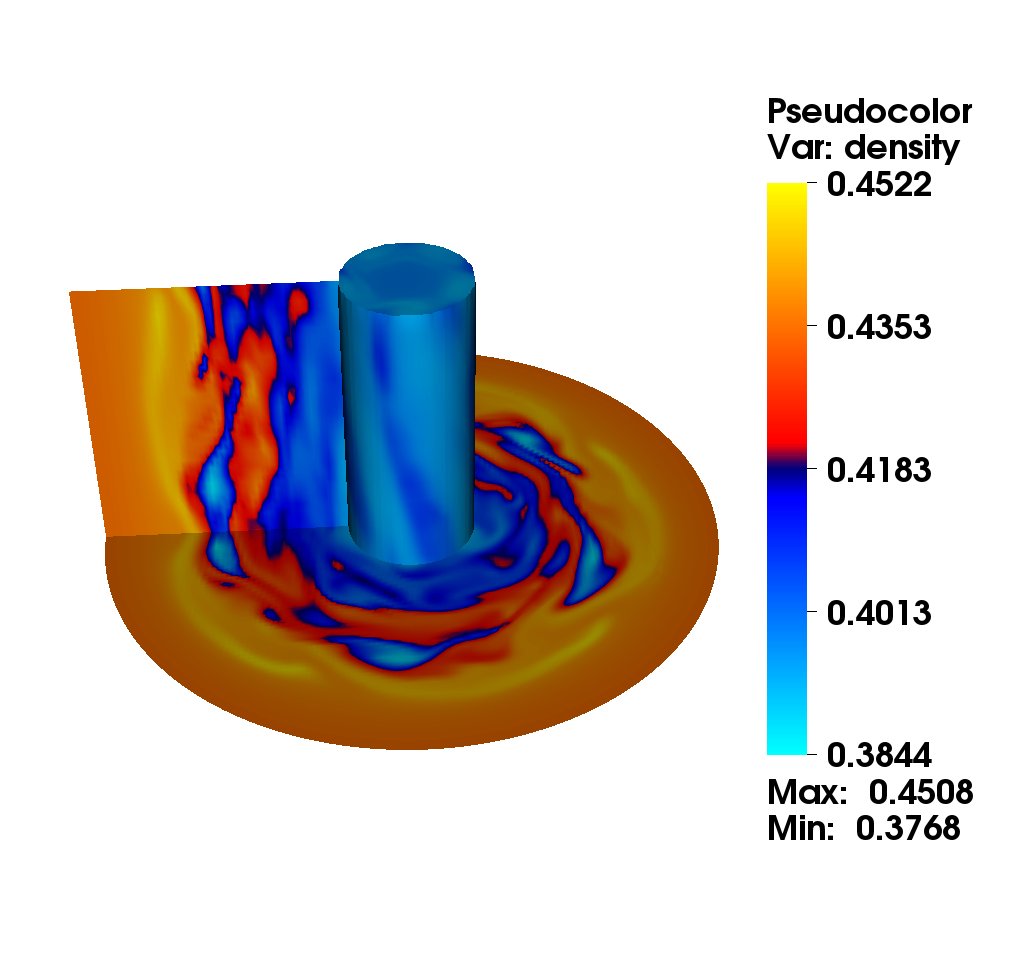} \\
  (c) $t=4000$                                   &    (d) $t=5000$  \\
  \includegraphics[width=6.5cm]{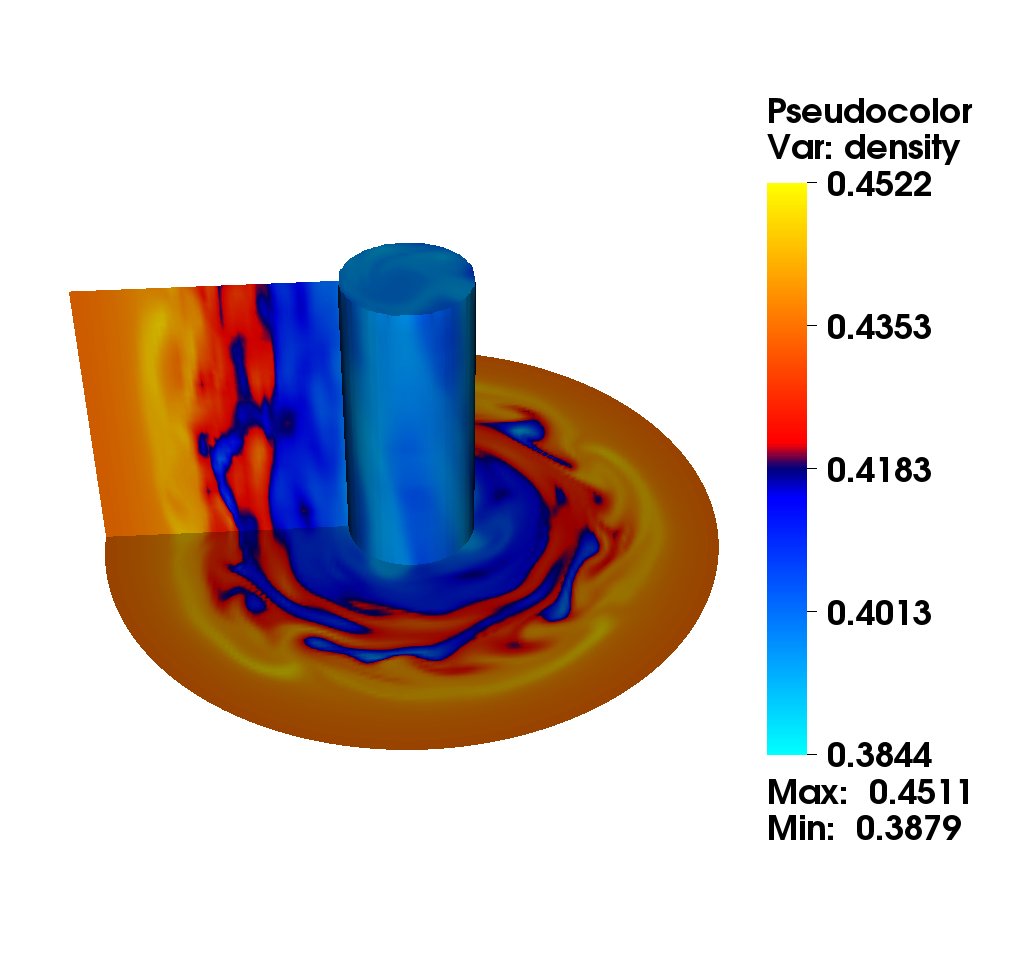} &    \includegraphics[width=6.5cm]{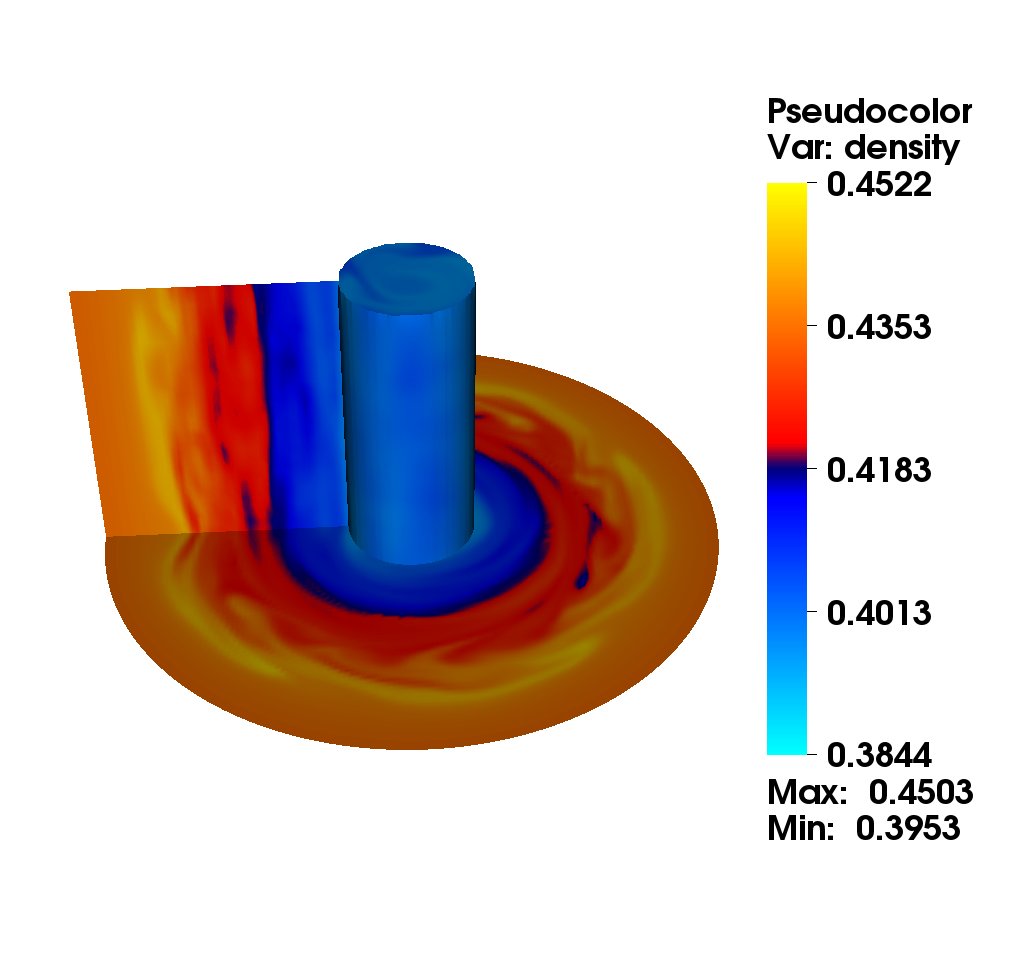} \\
  (e) $t=6000$                                   &       (f) $t=8000$
 \end{tabular}
\caption{\label{fig:DK_distribution_fun}Evolution of ion turbulence. The distribution function is shown for the velocity $v_\|=0$. The mesh size is $n_x=n_y=128,n_z=32,n_v=65$. Mixed  Semi-Lagrangian/finite difference method is used.}
 \end{center}
\end{figure}

In Figure~\ref{fig:DK_distribution_fun}, we show the evolution of distribution function $f$ at $v=0$. We see that till the time $t\leq2000$ the instability can not be identified very clearly, that is why the Semi-Lagrangian method can be applied in the linear phase. At time $t=3000$, we reaches the saturation point, and five vortices are developed. These vortices rotate and create small filaments. At time $t=4000$, more small structures appear in the distribution function. At this moment, the Semi-Lagrangian method can not conserve well the invariant quantities, while the conservative finite difference method performs much better. Finally, the instability continues to develop  small structures of the distribution function till the mesh size, and we attain to a relatively steady state.

\section{Conclusion and perspective}
\label{sec:conc}
\setcounter{equation}{0}

In this paper, we have presented an efficient algorithm for long term plasma simulations.
We first derive the 4D drift-kinetic and the 2D guiding-center models, and present their conservative properties. The Hermite WENO reconstructions are applied for solving the Vlasov equations, which was proved to be robust~\cite{YF15} in computational performance. Moreover, to adapt the arbitrary computational domain, we discretize the models on Cartesian meshes, and the special numerical methods for the boundary conditions, as the inverse Lax-Wendroff method for the Vlasov equation~\cite{bibFY1} and the extrapolation method for the poisson equation~\cite{bibFY2}, are proposed.

Next, we solve the guiding-center model on a D-shape domain. A steady state solution $(\phi_0,\bar\rho_0)$ is found numerically. Then we perturb the steady state density $\bar\rho_0$ along the streamline, and use this perturbed density $\bar\rho$ as the initial condition for the guiding-center model. We observe that the difference of density $\delta\rho=\bar\rho-\bar\rho_0$ revolves, and the filaments appear for long term simulation.

Finally, we simulate the 4D drift-kinetic model with the mixed methods, {\it i.e.} the semi-Lagrangian method in linear phase and finite difference method during the nonlinear phase. Numerical results show that the mixed method is efficient and accurate in linear phase and  it is much stable during the nonlinear phase. Moreover, it preserves well the conservative properties. We thus conclude that our mixed method is efficient for realistic and high dimensional plasma turbulence simulations.


\section*{Acknowledgment} 
 Both authors are  partially supported by the European Research Council
ERC Starting Grant 2009,  project 239983-\textit{NuSiKiMo} and the
Inria project Kaliffe. Chang YANG is also supported by National Natural Science Foundation of China (Grant No. 11401138).

\bibliographystyle{plain}

\end{document}